%% file: ARXIV_ex_article_new_resub.tex
\documentclass[review,onefignum,onetabnum,final]{siamonline250211}

\usepackage{lineno}
\input{ARXIV_ex_shared}

\ifpdf
\hypersetup{
 pdftitle={Efficient decomposition of Forman-Ricci curvature over Vietoris-Rips complexes and Data Applications
},
  pdfauthor={Danillo Barros de Souza,Jonatas Teodomiro, Fernando A. N. Santos, Mengjun Ding, Weiqiang Sun, Mathieu Desroches, J\"urgen Jost and Serafim Rodrigues}
}
\fi

\newtheorem{prop}{Proposition}
\usepackage{cleveref}
\usepackage{array}
\usepackage{booktabs} 
\usepackage{xcolor,hyperref}
\hypersetup{colorlinks,breaklinks,
            linkcolor=blue,urlcolor=blue,
            anchorcolor=blue,citecolor=blue}

\usepackage{todonotes}

\usepackage[normalem]{ulem}
\headers{F.R.C. in V.R. Complexes and data applications}{D.B.S. et al}
\begin{document}
\nolinenumbers
\maketitle

% REQUIRED
%\begin{abstract}
%Discrete Forman-Ricci curvature (FRC) has been shown as an efficient geometric hallmark for transitions in networks-based approaches, which can also be extended to higher dimensional computations in simplicial complexes. Aside from that, a natural computational complexity is carried along with computations, which are the biggest pitfalls for applications of geometric tools to distance-based network constructions. Furthermore, there is no local geometric descriptor that can be naturally extended to a global computation, as performed in the Gauss-Bonnet theorem. This issue motivated us to extend the derivations of set-theoretical approaches to design efficient algorithms for Vietoris-Rips (VR) complexes. As a result, it allowed as to compute higher-order FRC in function of the distance between simplices, which we refer as a data geometrization. Finally, we benchmark the data geometrization process in point-cloud data and breast cancer datasets and compare the results with the classical classification methods. Impressively, our approach uncovers crucial information from the data geometry that is not usually perceived by the current data classification methods. Our findings not only pave the way for geometric computations in VR complexes but also reveal a missing point for efficient data classifications: the intrinsic geometry.
%\end{abstract}
\begin{abstract}
%Discrete Forman-Ricci curvature (FRC) %has been shown as an efficient geometric hallmark for transitions in networks-based approaches, which can also be extended to higher-dimensional computations in simplicial complexes. In this work, we extend the computation of the FRC to local computations such as in the Gauss-Bonnet theorem, as well as the derivations of set-theoretical approaches to design an efficient algorithm for computing FRC in Vietoris-Rips (VR) complexes. {This allows us} to compute higher-order FRC {as a} function of the 
 %filtration distance in VR complexes, which we refer to as data geometrization. Finally, we benchmark the data geometrization process in {both} synthetic point-cloud data and breast cancer data. {What's more, we use} geometry as an input {in the context of} a {state-of-the-art} data classification method. Our approach {then} uncovers crucial information from the data geometry that is not usually perceived by the current data classification methods. Our findings pave the way for geometric computations in VR complexes and reveal a {key} missing {element in} efficient data classifications: {considering} the geometry behind the statistical information.
Discrete Forman-Ricci curvature (FRC) is an efficient tool that characterizes essential geometrical features and associated transitions of real-world networks, extending seamlessly to higher-dimensional computations in simplicial complexes. In this article, we provide two major advancements: First, we give a decomposition for FRC\textcolor{black}{that inhently  allows a} local computations of FRC. Second, we construct a set-theoretical proof enabling an efficient algorithm\textcolor{black}{s} for the local computation of FRC in Vietoris-Rips (VR) complexes. Our findings open new avenues for geometric computations in VR complexes and highlight an essential yet under-explored aspect of data \textcolor{black}{ analysis and visualisation}: the geometry underpinning statistical patterns.
\end{abstract}

% REQUIRED

% REQUIRED
\section{Introduction}
%Network analysis~\cite{newman2003structure} has been applied as an effective approach {to analyse complex} data, with a wide range of models~\cite{zhou2020graph} and applications to community detection~\cite{fortunato2010community}. These approaches are often used {to improve} current machine learning techniques~\cite{jin2021survey}. Topological data analysis \cite{zomorodian2012topological} has emerged in the last couple of decades as a powerful tool for enhancing the signal-to-noise of data features, and it has made crucial contributions across fields, in particular in neuroscience \cite{gunnarcarlsson2018towards}. The discrete Ricci curvatures~\cite{forman2003bochner,Lin2011RicciGraphs,samal2018comparative} have proven to be a set of powerful geometric descriptors of network clusters and they have been applied to numerous fields. To name a few: in the study of stock market fragility~\cite{sandhu2016ricci}, in neuroscience~\cite{chatterjee2021detecting}, epidemiology~\cite{de2021using}. In particular, the discrete Forman-Ricci curvature (FRC)~\cite{forman2003bochner} is a powerful geometric descriptor in complex networks. 
Network-based analysis \cite{newman2003structure} provides a versatile and effective framework for data mining \cite{zhou2020graph} and thus applies to a wide range of applications, particularly for data with complex relationships and hierarchies, such as community detection \cite{fortunato2010community}. Recently, it has been instrumental in advancing machine learning techniques \cite{jin2021survey}. Beyond the dyadic network framework, Topological data analysis (TDA) has recently emerged as the leading approach to studying higher-order relationships and hierarchies in data, as well as enhancing the signal-to-noise ratio of data features \cite{zomorodian2012topological}. As a consequence, TDA has permeated across several fields, including neuroscience \cite{gunnarcarlsson2018towards}. TDA has also advanced with the study of geometric properties of data (Topological and Geometrical Data Analysis - TGDA). In this context, discrete Ricci curvatures have been shown to be a powerful geometric descriptor of networks, enabling, for example, efficient detection of network clusters \cite{forman2003bochner,Lin2011RicciGraphs,samal2018comparative}. Its application is now widespread in several fields, such as stock market fragility \cite{sandhu2016ricci}, neuroscience \cite{chatterjee2021detecting}, and epidemiology \cite{de2021using}, to name a few. A step further is given by the discrete Forman-Ricci curvature (FRC), which provides a powerful geometric descriptor of complex networks\cite{forman2003bochner}. 

Despite the relevance of geometric descriptors, topological approaches, such as those based on persistent homology~\cite{edelsbrunner2008persistent,zomorodian2005topology} prioritize topological descriptors of data while largely disregarding the geometric perspective. This is partly explained due to the computational complexity of geometric descriptors, particularly in the context of constructing higher-order networks.

To overcome these challenges, we leverage set theory~\cite{enderton1977elements,jech2003set,pawlak2002rough,zimmermann1985applications}, which has proven effective in optimizing discrete computational algorithmic methods. For instance, our recent works unveiled the construction of efficient set-theoretical approaches for computing geometric invariants~\cite{de2023efficient,desouza2025alternativesettheoreticalalgorithmsefficient}. Building on these results, we explore geometric computations in greater detail and identify computational patterns in the construction of Vietoris-Rips (VR) complexes that reduce the complexity of the FRC computation. We find a decomposition strategy that enables the local computations of FRC from local network neighbourhoods. This permits the geometric local update to be possible with minor numerical increments and efficient computational complexity.

 To clarify the rationale behind our approach, we give an overview of our computational strategy. In traditional homology barcode computations, the topological invariant is recomputed as a function of the cutoff distances of VR complexes, once the global structure needs to be recomputed from scratch. For example, in order to recompute the Betti numbers~\cite{bochner1949curvature}, the clique's neighbourhood needs to be recomputed, as well as the boundary operators. From a geometric point of view, the local neighbourhood changes as the new structures are added as a function of the distance, and may imply re-computations of local geometry as well. To reduce the computational complexity, we develop an alternative set-theoretical approach that updates the numerical computation of these curvatures instead of recomputing the FRC for the updated network. \textcolor{black}{Also, TDA has proven to be an effective framework for analyzing high-dimensional and noisy datasets, owing in part to the stability of topological descriptors such as persistence diagrams under small perturbations of the input data~\cite{cohen2005stability}. Motivated by this setting, we develop an efficient computational framework for Forman--Ricci curvature on Vietoris--Rips complexes.} As a proof-of-concept, we test our novel computational method on various datasets with the aim of exploring and understanding the effect of geometric approaches in VR complexes. This also includes noise sensitivity studies \textcolor{black}{in synthetic data}. Here, we denote our methodological procedure \textit{data geometrization}. Our studies validate our approach and, crucially, introduce a novel descriptor suited for high-dimensional datasets.

The rest of this article is organised as follows. In~\Cref{sec:FRC_in_VR}, we give a brief review of the theoretical underpinnings of networks and VR complexes, define the FRC and examine associated optimizations methods in past works. Then in~\Cref{sec:results}, we benchmark our novel algorithm to compute FRC in synthetic \textcolor{black}{ datasets, and further compare with the Ollivier-Ricci curvature outputs } in real datasets. Finally, in~\Cref{sec:conclusion}, we discuss the details of our findings and give a few perspectives on this topic.
%\MDcom[inline]{No need for discussion and conclusion, I would merge them, especially that they are both rather short.}

%-----------------------------------------
\section{Theoretical background}
\label{sec:FRC_in_VR}
%-----------------------------------------
In this section, we provide a brief background of graphs and the VR complexes. Subsequently, we define the discrete FRC curvature. We refer the reader to classic literature for more details on these concepts~\cite{edelsbrunner2022computational,forman2003bochner,zomorodian2010fast,zomorodian2005topology}.
%-------------
\subsection{Vietoris-Rips complexes from networks}
\label{sec:graphs_and_VR}
%-------------
A simple weighted undirected graph is a pair $G=(V,E)$, together with a weight function $w:E\rightarrow\mathbb{R}$, where $V:=V(G)$ is the (finite) set of nodes of $G$ and $E:=E(G)$ is the set of edges connecting nodes in $G$, such that the following equation is satisfied:
\begin{eqnarray}
    E\subseteq\{\{x,y\}|\, x,y\in V,\, x\neq y\}.
\end{eqnarray}
The neighbours of a node $x \in V$ will be denoted by $\pi_x$ and are defined by
\begin{eqnarray}
\label{eq:node_neighbourhood}
\pi_x=\{y \in V|\, \{x,y\}  \in E\}.
\end{eqnarray}
Similarly, the neighbours of a set of nodes $\alpha=\{x_1,\hdots,x_n\}$ are defined by
\begin{eqnarray}
\label{eq:cell_neighbourhood}
\pi_\alpha=\bigcap_{x\in \alpha}\pi_x.
\end{eqnarray}
A {\textit{clique complex}} is the set of all complete subgraphs of $G$. We say that these subgraphs are its simplices, and, when such a subgraph has $d+1$ nodes, we call it a $d$-face. We define the set of all $d$-faces by $C_d$. We say that a $d$-face has dimension $d$ and that the dimension of the clique complex is the highest dimension among its faces. We can define this clique complex as the union $C=\bigcup_{d} C_d$. {We see that $V(G)$, together with C, is an \textit{abstract simplicial complex}, i.e., they satisfy the following conditions:}
\begin{enumerate}
    \item For each $v \in {V(G)}$, $\{v\} \in {C};$
    \item If $\gamma \subseteq {\alpha}$ and  ${\alpha}\in {C}$, then $\gamma \in {C}$.
\end{enumerate}
We define the  \textit{Vietoris-Rips complex} constructed from $C$ as a function of the radius distance $\varepsilon$ by
\begin{eqnarray}
    \label{eq:VR_complex}
    \text{VR}_{C}(\epsilon):=\{\sigma \subseteq C\,|\, \text{diam}\, \sigma \leq \varepsilon \},
\end{eqnarray}
where $\text{diam }\sigma = \max \{w_e\,|\, e \in \sigma \}$, $w_e=w(e)$, and $\varepsilon$ is the \textit{cutoff} distance for generating the simplicial complex. Finally, let $\varepsilon_1,\varepsilon_2,\cdots , \varepsilon_k$ be the diameters of the simplices in the simplicial complex such that $\varepsilon_1\leq\varepsilon_2 \leq \cdots \leq \varepsilon_k$. Then it is clear that
\begin{eqnarray}
    \label{eq:sequence_VR_classic}
    \emptyset \subseteq \text{VR}(\varepsilon_1) \subseteq \cdots \subseteq \text{VR}(\varepsilon_k)=C.
\end{eqnarray}
\autoref{fig:FRC_filtration_process} (top) elucidates the construction of a VR complex with an intuitive example.
Next, we define the FRC and summarise the key results from our previous work, which will serve as the foundation for proposing a novel, efficient formulation.

%-------------
\subsection{Discrete Forman-Ricci curvature}
\label{sec:FRC}
%-------------
The \textit{boundary} of a $d$-face $\alpha$ is the set of all {$(d-1)$-faces} in $\alpha$, and it is denoted by $\partial(\alpha)$. When $\gamma \in C_{d-1}$ is contained in the boundary of $\alpha$, we denote it by $\gamma < \alpha$; alternatively, $\alpha$ containing $\gamma$ in its boundary will be denoted by $\alpha > \gamma$.
Two faces $\alpha_1,\alpha_2 \in C_d$ are said to be \textit{neighbours} if at least one of the following conditions is satisfied:
\begin{enumerate}
    \item there exists a ($d-1$)-face $\gamma$ such that $\gamma<\alpha_1,\alpha_2$,
     \item there exists a ($d+1$)-face $\beta$ such that $\alpha_1,\alpha_2<\beta$,
  \end{enumerate}
and we define by $N_\alpha$ the set of all neighbours of $\alpha$.
  \textcolor{black}{For simplicial complexes, }we say that $\alpha_1$ and $\alpha_2$ are \textit{parallel neighbours} if conditions $1.$ is satisfied, but not $2.$.  If both conditions are satisfied simultaneously, then $\alpha_1$ and $\alpha_2$ are said \textit{transverse neighbours}.

The original formulation for higher-order FRC, defined for CW-complexes, originates from \cite{forman2003bochner}. In the special case of simplicial complexes, we can express it as proposed in~\cite{de2023efficient} as follows
%The original formulation for higher-order FRC was derived in \cite{forman2003bochner} and it is given by
%
\begin{eqnarray}
\label{eq:original_FRC}
    \mathrm{F}(\alpha)=|H_{\alpha}| + (d+1) - |P_\alpha|,
    \end{eqnarray}
where $H_\alpha$  is the set of $(d+1)$-faces containing $\alpha$ in its boundary, and $P_\alpha$ is the set of parallel  neighbours of $\alpha$.
We therefore define the $d$-th Forman-Ricci curvature (or the $d$-FRC) of a non-empty simplicial complex as follows: 
\begin{eqnarray}\label{eq:avg_FRC}
    F_d(C)=\dfrac{1}{|C_d|}\sum_{\alpha\in C_d}F(\alpha),
\end{eqnarray}
assuming that $C_d \neq \emptyset$.
%Formally, in the case of degenerate simplicial complexes, we do not define curvature.
For convenience of numerical computations, we define $F_d(\emptyset):=0$. This will facilitate the formalism of our approach in this article. Finally, we define $F(\alpha)$ using a key result from our previous work~\cite{de2023efficient} as follows:
\begin{eqnarray}\label{eq:new_formula_FRC}
    \mathrm{F}(\alpha)=(d+2)\cdot\left|\bigcap_{\gamma \in \partial \alpha}\pi_\gamma\right|+2\cdot (d+1)-\sum_{\gamma \in \partial \alpha}|\pi_\gamma|,
\end{eqnarray}
where $\pi_\gamma$ is defined as in~\eqref{eq:cell_neighbourhood}.
%
%--------------------
\section{Results}
\label{sec:results}
%--------------------
We subsequently outline our main theoretical enhancements alongside their computational implications. Finally, to demonstrate its validity, we apply it to synthetic data.
%-------------
\subsection{Alternative formulation for local FRC computations}
\label{sec:local_FRC}
%-------------
 The classic FRC formulation  is originally defined from local neighbourhood interactions as in~\eqref{eq:original_FRC} and~\eqref{eq:avg_FRC} ( which uses the average of local $FRC$ of neighbouring nodes), and is already able  to recover the geometry from local to global scale \cite{bloch2014combinatorial,weber2018coarse}. This is usually possible in topological approaches, e.g., the Euler characteristic computation obtained from the Knill curvature~\cite{knill2020index}. Taking this facts into account, we investigated alternative ways to consider the decomposition of global curvature information to local computations . To this end, we constructed an alternative local computation of FRC inspired by the Gauss-Bonnet theorem for simplicial complexes \cite{knill2012discrete,knill2024gauss}. More precisely, we define the \textit{local Forman-Ricci curvature} to the nodes $x\in V$ (for $x$ in $d$-faces) as follows:
\begin{equation}
    \label{eq:local_total_forman}
    \mathrm{f}_d(x)=\frac{1}{(d+1)\cdot |C_d|}\sum_{\alpha \supset x} \mathrm{F}_d(\alpha).
\end{equation}
Thus, the global FRC can be recovered from the local computations from~\eqref{eq:local_total_forman}, specifically, by the formula
\begin{eqnarray}\label{eq:geometric_gauss_bonnet}
    \label{eq:average_frc}
    \mathrm{F}_d(C)=\sum_{x\in V} \mathrm{f}_d(x).
\end{eqnarray}
The derivation of~\eqref{eq:geometric_gauss_bonnet} is detailed in~\ref{sec:demonstrations}. \autoref{fig:example_adding_face} elucidates the computation of FRC. \textcolor{black}{Such derivation provides an exact computation of the classic global FRC (left side) from the defined local FRC (right side -- derived from equation \eqref{eq:local_total_forman}) in function of the VR complex.}

%-------------
\subsection{A novel algorithm to efficiently compute FRC in a VR complex}\label{sec:novel_algorithm}
\begin{figure}
    \centering
    \includegraphics[width=\linewidth]{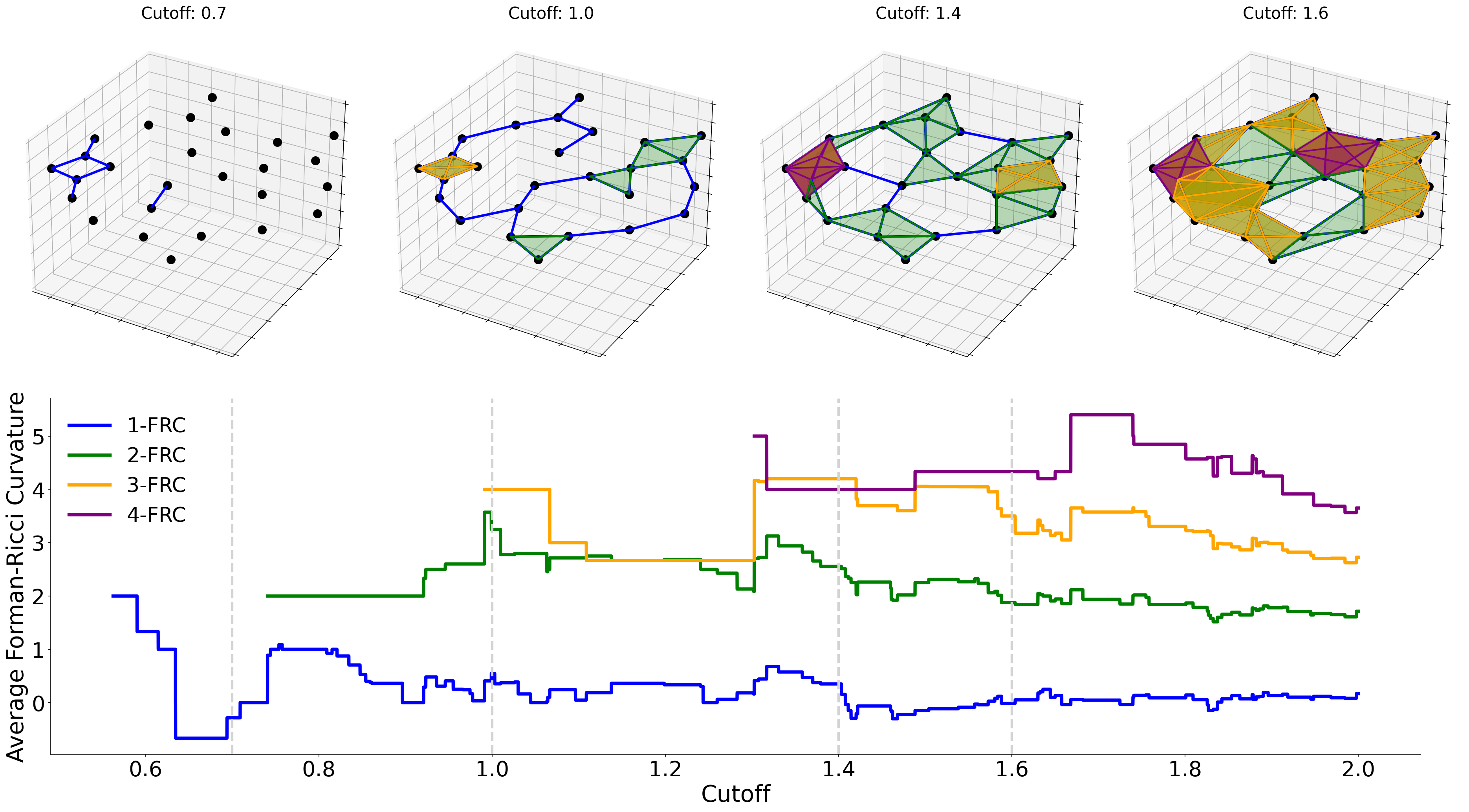}
    \caption{Evolution of a Vietoris--Rips (VR) filtration and the corresponding average Forman--Ricci curvature (FRC). Top: Sequence of VR complexes constructed from a three-dimensional point cloud at cutoff distances $\varepsilon=0.7$, $1.0$, $1.4$, and $1.6$. Increasing the cutoff distance progressively introduces additional simplices, thereby generating a nested filtration of simplicial complexes. Bottom: Average Forman--Ricci curvature as a function of the cutoff distance for simplices of dimensions $d=1$ (blue), $d=2$ (green), $d=3$ (yellow), and $d=4$ (red). The vertical dashed grey lines denote the cutoff values illustrated in the top panel, facilitating the correspondence between the geometric evolution of the VR filtration and the associated curvature values.}
    \label{fig:FRC_filtration_process}
\end{figure}

%-------------
Computations of higher dimensional geometric invariants in complex networks are not new. However, the computations on VR complexes are challenging, since the additional computation of several sequential networks may be necessary. \cref{fig:FRC_filtration_process} provides an example of FRC computation in a VR complex. Subsequently, we define such computations and explicitly point out the challenges that they entail. Then, we will propose a novel set-theoretical approach for efficiently computing high dimensional FRC in VR complexes filtrations (i.e., the FRC computation on high dimensional simplices).
Let $C_d=\{\alpha_1,\cdots,\alpha_k\}$ be a finite set (or sequence) of $d$-faces, with a positive weight function $\omega:C_d \rightarrow \mathbb{R} $ such that $\omega(\alpha_k)=w_{\alpha_k}$. We also assume that faces' weights are sorted in ascending order (\textit{i.e.,} $w_{\alpha_1}\leq w_{\alpha_2}, \leq \dots \leq w_{\alpha_k}$). We define $C^i_d=C_d^{i-1}\sqcup \{\alpha_i\}$, for $i>1$, with $C^0_d=\emptyset$ and $C^{1}_d=\{\alpha_1\}$. It is clear that the sequence of simplicial complexes $(C^i_d)_i$ are such that
\begin{eqnarray}
\label{eq:sequence_of_simplicial_complexes}
\emptyset \subseteq C^1_d\subseteq C^2_d\subseteq \cdots \subseteq C^k_d=C_d.
\end{eqnarray}
Therefore, \eqref{eq:sequence_of_simplicial_complexes} defines a VR complex.
A natural and straightforward FRC computation in a VR complex can be accomplished by simply computing the sequence
\begin{eqnarray}
\label{eq:sequence_of_FRC_computations}
0=:\mathrm{F}(\emptyset), \,\mathrm{F}(C^1_d),\, \mathrm{F}(C^2_d), \, \cdots , \,\mathrm{F}(C^k_d)=\mathrm{F}(C_d).
\end{eqnarray}
However, this would imply an algorithm that requires exhaustive recomputation of the FRC for each simplicial complex $C^i_d$.
Instead, we find that linking equations \eqref{eq:sequence_of_simplicial_complexes} and \eqref{eq:sequence_VR_classic} enables the construction of an optimal algorithm for computing the FRC of higher-order faces across a filtration. This observation leads us to propose an optimal algorithm which takes as input the sorted simplicial complex $C$ and the maximum face dimension $d_{max}$. Such an algorithm will be detailed in this section.

Henceforth, for all $\alpha \in C^i_d$, we will denote $N_\alpha^{i}:=N_\alpha(C^i_d),$ $T_\alpha^{i}:=T_\alpha(C^i_d),$ $P_\alpha^{i}:=P_\alpha(C^i_d)$ and $H_\alpha^{i}:=H_\alpha(C^i_d)$ the set of neighbours, transverse neighbours, parallel neighbours and higher-order faces containing $\alpha$ with regards to $C^i_d$, respectively. In the special case where $\alpha=\alpha_i$, we will denote $N^i_{\alpha_i}$ by $N^i$ and $P^i_{\alpha_i}$ by $P^i$. We also define $\mathrm{F}^i(\alpha)$ by the computation of FRC to $\alpha$ at the $i$-th step, i.e., when $\alpha\in C^i_d,\;\alpha \neq \alpha_i$. More precisely, we define: $\mathrm{F^{i}(\alpha)}=|H^i_\alpha|+(d+1)-|P^i_{\alpha}|, \forall \alpha \in C^{i-1}_d$.

From Proposition~\ref{prop:before_after} (see section~\autoref{sec:demonstrations}), if $\alpha^i$ is the new neighbour of some $\alpha \in C^i_d$, then it can be either a parallel or transverse neighbour of $\alpha$. If $\alpha_i$ is a new parallel neighbour of  $\alpha$, then we have $\mathrm{F}^{i}(\alpha)=\mathrm{F}^{i-1}(\alpha)-1$. Otherwise, we have $\alpha$ as a new transverse neighbour and thus $\mathrm{F}^{i}(\alpha)=\mathrm{F}^{i-1}(\alpha)+(d+1)$. For illustration, \autoref{fig:example_adding_face} provides an example for edges ($d=1$) and triangles ($d=2$). 

This finding, together with the new formulation provided in~\eqref{eq:new_formula_FRC} from~\cite{de2023efficient} motivated us to derive and design an algorithm for computing FRC in a filtration of a VR complex.
In order to facilitate the implementation of our algorithm, we define $\Delta\mathrm{F}(\alpha):=\mathrm{F}^i(\alpha)-\mathrm{F}^{i-1}(\alpha)$ as an auxiliary function to be applied over the neighbours $\alpha \in N_{\alpha_i}^i$, which can be re-written as follows:
%\textcolor{red}{Need to improve for neighbours}
%
\begin{eqnarray}
    \label{eq:delta_function}
    \Delta\mathrm{F}(\alpha)=
    \left\{
    \begin{array}{ll}
-1, &  \alpha \in P^i \\
(d+1), &  \text{otherwise}
\end{array}
\right.
.
\end{eqnarray}
Using our results from~\cite{de2021using} (see also~\autoref{sec:previous_results}), this formula induces an equivalent set-theoretical formulation as a function of the nodes $x \in \partial(\alpha)$ as follows:
\begin{eqnarray}
    \label{eq:delta_function_node}
    \delta(x)=
    \left\{
    \begin{array}{ll}
-1, &  x \notin \pi_{\alpha_i} \\
(d+1), &  \text{otherwise}
\end{array}
\right.
.
\end{eqnarray}
The core idea behind formulas~\eqref{eq:delta_function} and~\eqref{eq:delta_function_node} is to examine the neighbourhood of the new face at step $i$ and locally iterate through the neighbouring faces to assess the contribution of the new face added to the curvature. This is achieved simply by comparing the curvature before and after adding the new face. In particular, equation~\eqref{eq:delta_function_node} allows this verification by identifying the neighbourhood of the current face through a set-theoretical representation, as performed in our previous work in \cite{de2023efficient}. The above formulation, particularly through~\eqref{eq:delta_function_node}, enables the implementation of an efficient algorithm for computing the FRC on VR complexes, using a decision tree based on the neighbourhood of each newly added face $\alpha_i$. This leads us to our algorithm \eqref{alg:FRC_raw} (see~\autoref{sec:algorithms} for more alternative algorithms).

\begin{algorithm}
\caption{Compute Local and Global Forman-Ricci Curvature (FRC)}
\label{alg:FRC}
\begin{algorithmic}[1]
\STATE \textbf{Input:} Distance matrix $D$, maximum dimension $d_{\max}$, maximum distance $\mathrm{max\_dist}$, precision $p$
\STATE \textbf{Output:} Average local FRC and total FRC for nodes $V$
\STATE

% Step 1: Initialize parameters
\STATE Initialize node set $V$ with $|V|$ nodes (if not provided, assume $V = \{0,1, \ldots, |D|\}$)
\STATE Set cutoff step size $\delta \leftarrow 10^{-p}$ and compute $\mathrm{cutoffs} \leftarrow [0, \mathrm{max\_dist}]$ with step $\delta$
\STATE Initialize: 
\STATE $\mathrm{C}_d \leftarrow 0$, $\mathrm{F}_d \leftarrow 0$ (global curvature), $\mathrm{nF}_d(x) \leftarrow 0 \, \forall x \in V$
\STATE $\mathrm{Neighborhood} \, N_d \leftarrow \{\}$ for $d \in \{1, \dots, d_{\max}\}$

% Step 2: Generate cliques and compute FRC
\STATE \textbf{if} $d_{\max} = 1$:
    \STATE \quad Generate edge list $E$ and corresponding distances $w_{ij}$ sorted by increasing order
\STATE \textbf{else:}
    \STATE \quad Use \texttt{cliques\_gudhi} to generate cliques $\alpha$ with $w_{\alpha} \leq \mathrm{max\_dist}$
\STATE

\FOR{each clique $\alpha$ with weight $w_\alpha$}
    \STATE $d \leftarrow |\alpha| - 1$ \COMMENT{Dimension of the clique}
    \STATE $\mathrm{C}_d \leftarrow \mathrm{C}_d + 1$ \COMMENT{Count the clique}
    
    % Step 3: Update boundary neighborhoods
    \STATE Compute boundary $\partial(\alpha) \leftarrow \{B \,|\, B \subset \alpha, |B| = d\}$
    \FOR{each boundary $B \in \partial(\alpha)$}
        \STATE Update $\mathrm{Neighborhood}$: $N_d[B] \leftarrow N_d[B] \cup (\alpha \setminus B)$
    \ENDFOR
    
    % Step 4: Compute FRC for the clique
    \STATE Compute $H \leftarrow \bigcap_{B \in \partial(\alpha)} N_d[B]$ \COMMENT{Intersection of neighbors of boundaries}
    \STATE Compute total neighbors $n \leftarrow \sum_{B \in \partial(\alpha)} |N_d[B]|$
    \STATE Compute curvature: $\mathrm{f}_d \leftarrow (d+2) \cdot |H| + 2 \cdot (d+1) - n$

    % Step 5: Update local FRC for nodes
    \FOR{each node $x \in \alpha$}
        \STATE $\mathrm{nF}_d(x) \leftarrow \mathrm{nF}_d(x) + \mathrm{f}_d / (d+1)$
    \ENDFOR

    % Step 6: Update FRC for boundary nodes and neighborhoods
    \FOR{each boundary $B \in \partial(\alpha)$}
        \FOR{each node $x \in N_d[B] \setminus \alpha$}
            \STATE Compute $\delta \leftarrow \delta(x \in H, d)$ \COMMENT{Transverse or parallel neighbor}
            \FOR{each node $y \in B$}
                \STATE $\mathrm{nF}_d(y) \leftarrow \mathrm{nF}_d(y) + \delta / (d+1)$
            \ENDFOR
            \STATE $\mathrm{nF}_d(x) \leftarrow \mathrm{nF}_d(x) + \delta / (d+1)$
            \STATE Update total curvature: $\mathrm{f}_d \leftarrow \mathrm{f}_d + \delta$
        \ENDFOR
    \ENDFOR

    \STATE Update global curvature: $\mathrm{F}_d \leftarrow \mathrm{F}_d + \mathrm{f}_d$

    % Step 7: Print or store results when cutoff changes
    \IF{cutoff $w_\alpha$ changes from previous $w$}
        \IF{$\mathrm{C}_d > 0$}
            \STATE Output: $(w, \mathrm{C}_d / \mathrm{combs}_d, \mathrm{F}_d / \mathrm{C}_d, \{\mathrm{nF}_d(x) / \mathrm{C}_d, \forall x \in V\})$
        \ELSE
            \STATE Output: $(w, 0, 0, \{0, \forall x \in V\})$
        \ENDIF
        \STATE Update $w \leftarrow w_\alpha$
    \ENDIF
\ENDFOR

% Step 8: Finalize results
\STATE Fill in missing cutoff results by interpolating last valid curvature values
\STATE Return final FRC values: Average local FRC and total FRC for nodes
\end{algorithmic}
\end{algorithm}
This novel derivation enables an efficient and dynamic computation of the FRC as a function of the cutoff distance by computing updates of local curvature on the fly. Furthermore, such computation recovers the same classic global FRC information, once it uses the formulation provided in equation~\eqref{eq:geometric_gauss_bonnet}. 
We additionally estimate the time  and space complexity of the algorithm \eqref{alg:FRC}.
The time complexity of the algorithm is dominated by the nested loops used to update the curvatures of the faces in each cell. For a $d$-cell $\alpha$, these loops require $O(d^2\mathcal{B}_\alpha)$ operations, where $\mathcal{B}_\alpha=\max\{|N_d[B]|:B\in\partial(\alpha)\}$ and $N_d[B]$ denotes the set of nodes stored in the $d$-dimensional neighbourhood of the boundary face $B$ in Algorithm~\eqref{alg:FRC}. Therefore, processing all cells in $C$ up to dimension $d_{\max}$ has time complexity $O(|C|\cdot d_{\max}^2\cdot \Delta(G))$, where $\Delta(G)$ is the maximum degree of the graph $G$. The space complexity is dominated by the neighborhood of faces, with $O(|C|\cdot d^2_{\max})$ and the output stored, with $O(\lceil m_d10^\delta\rceil nd_{\max})$, where $m_d$ is the maximum distance, $n$ is the number of vertices, and $\delta$ is the step for incrementing the cutoff values in the VR filtration. Therefore, the total space complexity could reach $O(d_{\max}(|C|\cdot d_{\max}+\lceil m_d \cdot 10^{\delta}\rceil  n))$.
 The current implementation is integrated to FastForman package \cite{Barros_de_Souza_FastForman_-_An_2024}. Finally, we provide with an example a Python implementation of our algorithm that benchmarks the methods used for computing FRC from the local derivations defined in our work (see \cite{souza2025fastforman_benchmark}). 
\subsection{Application of Data Geometrization}
\label{sec:results_data}
\begin{table}[]
\label{tab:classic_table}
\centering
\begin{tabular}{|>{\centering}p{25mm}|>{\centering}p{25mm}|>{\centering}p{25mm}|p{5mm}|>{\centering\arraybackslash}p{25mm}|}
\hline
Observation & Feature 1 & Feature 2 & $\dots$ & Feature n \\ 
\hline
$x_1$ & $x_{11}$ & $x_{12}$ & $\dots$ & $x_{1n}$ \\ 
\hline
$x_2$ & $x_{21}$ & $x_{22}$ & $\dots$ & $x_{2n}$ \\ 
\hline
$\dots$ & $\dots$ & $\dots$ & $\dots$ & $\dots$ \\ 
\hline
$x_m$ & $x_{m1}$ & $x_{m2}$ & $\dots$ & $x_{mn}$ \\ 
\hline
\end{tabular}
    \caption{Classic Tabular Input Data: The input table is given as an $m \times n$ matrix (\textit{i.e.}, $m$ points in a $n$-dimensional space). More precisely, each observation point is in the shape $x_i=(x_{i1},x_{i2},\hdots,x_{in})$ for $i \in \{1,2,\hdots,m\}$.}
\end{table}
\begin{table}[!]
\label{tab:Geometrized_table}
\centering
\begin{tabular}{|>{\centering}p{25mm}|>{\centering}p{25mm}|>{\centering}p{25mm}|p{5mm}|>{\centering\arraybackslash}p{25mm}|}
\hline
Observation & $\varepsilon_1$ & $\varepsilon_2$ & $\dots$ & $\varepsilon_k$ \\ 
\hline
$\mathrm{f}_d(x_1)$ & $\mathrm{f}_d^{11}$ & $\mathrm{f}_d^{12}$ & $\dots$ & $\mathrm{f}_d^{1k}$ \\ 
\hline
$\mathrm{f}_d(x_2)$ & $\mathrm{f}_d^{21}$ & $\mathrm{f}_d^{22}$ & $\dots$ & $\mathrm{f}_d^{2k}$ \\ 
\hline
$\dots$ & $\dots$ & $\dots$ & $\dots$ & $\dots$ \\ 
\hline
$\mathrm{f}_d(x_m)$ & $\mathrm{f}_d^{m1}$ & $\mathrm{f}_d^{m2}$ & $\dots$ & $\mathrm{f}_d^{mk}$ \\ 
\hline
\end{tabular}
\caption{Geometrized Tabular Data: The new input table is provided from the geometric information (FRC) computed on each observation originally provided by \cref{tab:classic_table}. More precisely, each geometrized observation is a vector $\mathrm{f}_d(x_i)=(\mathrm{f}_d^{i1},\mathrm{f}_d^{i2},\dots,\mathrm{f}_d^{ik})$, for all $i\in \{1,2,\dots,k\}$, where $\mathrm{f}_d^{ij}$ is the computation of the $d$-th local FRC for the node $x_i$ restricted to $\text{VR}(\varepsilon_j)$. The global FRC is recovered from \eqref{eq:geometric_gauss_bonnet}, and the geometric output is similar to the content in \cref{fig:FRC_filtration_process}.}
\end{table}
%------------
Data science approaches aim at determining descriptive statistical summaries that synthesizes data, enabling inferences and predictions. Typically, this data is provided as in \autoref{tab:classic_table}. In our approach, we perform the FRC per observation (i.e., points in the aforementioned assumed coordinate system defined by the data features) by using equation~\eqref{eq:local_total_forman} as a function of each distance $\varepsilon$ (in this case, the Euclidean distance), which leads to \cref{tab:Geometrized_table}. The data geometrization concept has been applied to manyfold learning approaches, as in \cite{xu2022amcad}.
To validate our proposed approach, we consider both synthetic point cloud data and breast cancer data from public repositories~\cite{UCI_Wisconsin_Breast_Cancer,METABRIC_Kaggle}. Our algorithm and all data processing were implemented in the Python language~\cite{python}. To facilitate understanding, we also provide the pseudo-code to compute both local and global FRC in \autoref{sec:algorithms}. To generate the faces of the simplicial complexes, we use the package Gudhi~\cite{gudhi:2014,gudhi}. Following our approach, the computed curvature values (as a function of the cutoff distance) can then be used as features for revealing patterns and data visualisation. In particular, we employ the Uniform Manifold Approximation and Projection algorithm (UMAP \cite{umap-learn,umap}) to generate embeddings from the geometric outputs.

%------------------
\subsubsection{FRC computations on Synthetic point cloud data}
%------------------
We examine two sub-categories of synthetic point clouds, namely Random Geometric Graphs and Replicas of the Datasaurus datasets. The corresponding generated datasets are publicly available via \cite{FRC_VR_kaggle}. These two sub-cases are detailed below:

\begin{enumerate}
\item Random geometric graphs \cite{dall2002random,penrose2003random} with $n=100$ nodes and different box dimensions ($\text{dim} \in \{2,4,6,8,10,12,14,16,18,20\}$) and network densities ($\rho\in \{0.1,0.25,05,0.75,0.9\}$). To generate the random geometric graphs, we used the package NetworkX~\cite{networkx}. Noteworthy, we tested the sensitivity of the box dimension across the network density;
\item Replicas of Datasaurus datasets, (a synthetic point-cloud data repository) that generates different geometries with the same basic statistics~\cite{datasaurus}(see~\autoref{fig:datasaurus_dataset}). This is inspired by an alternative randomization algorithm that shuffles the points coordinates while preserving the average, the standard deviation and Pearson correlation of points coordinates~\cite{matejka2017same}. For more details, see~\autoref{sec:algorithms}. As above, we tested the sensitivity of the data recognition across datasets. 
\end{enumerate}
%
%
%\MDcom[inline]{This should be indicated at the end, in ``Data availability''}
%
\begin{figure}
    \centering
    \includegraphics[width=\linewidth]{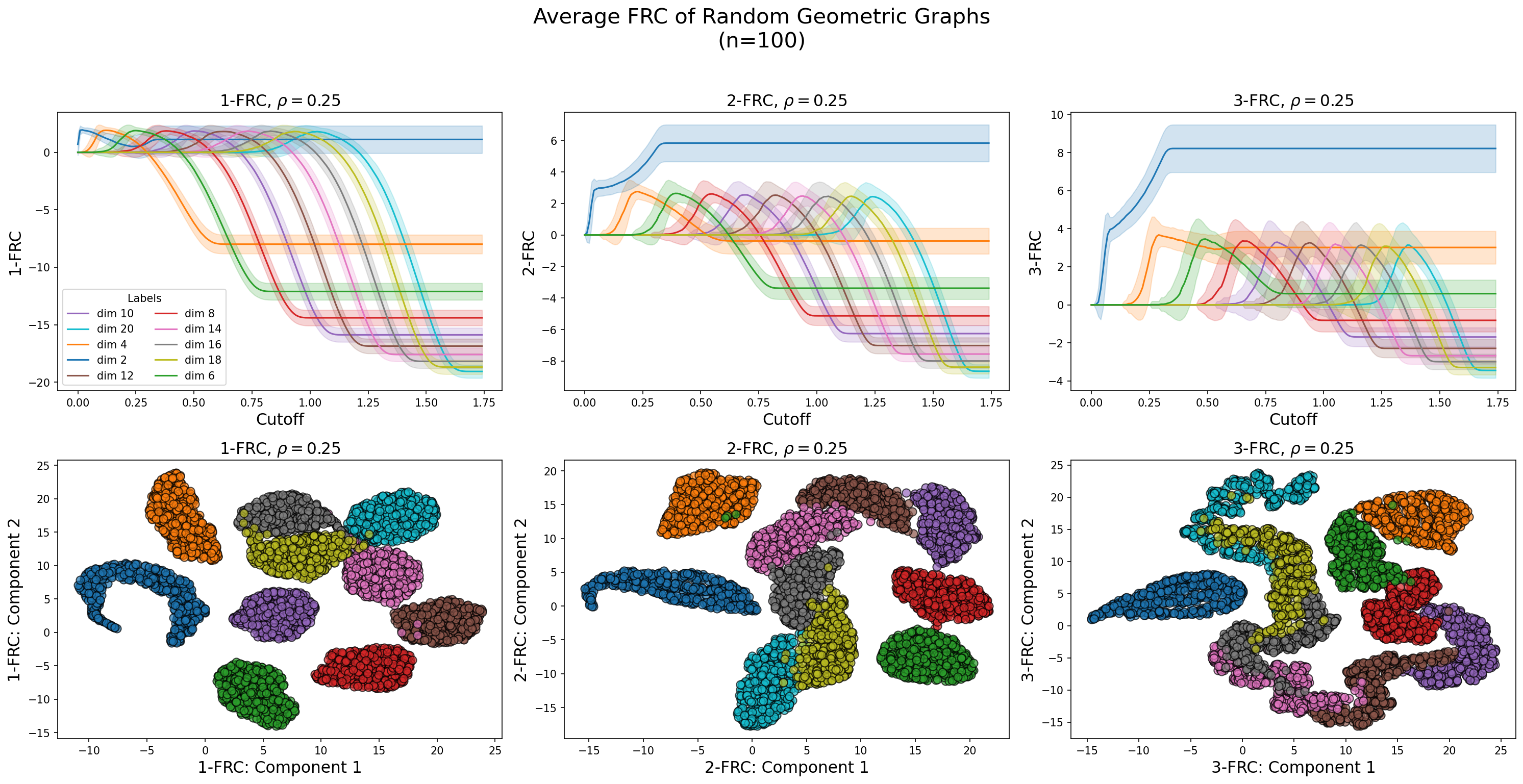}
    \caption{\textcolor{black}{Average Forman--Ricci curvature and the corresponding UMAP embeddings for random geometric graphs with varying box dimensions. Top: Average $d$-dimensional Forman--Ricci curvature ($d$-FRC) computed as a function of the filtration parameter for $d=1$ (first column), $d=2$ (second column), and $d=3$ (third column). Each curve corresponds to random geometric graphs generated with fixed edge density $\rho=0.25$ and different box dimensions (see legends). Bottom: Two-dimensional UMAP embeddings of the average $d$-FRC curves shown in the corresponding top panels. The embeddings were computed using the Euclidean metric with the maximum permitted minimum distance parameter ($\mathrm{min_dist}=1.0$), thereby maximizing the separation between embedded points to assess the robustness of the average $d$-FRC descriptors to noise.}}
    \label{fig:UMAP_and_FRC_RGG_den_025}
\end{figure}
Starting with the case of random geometric graph, in~\cref{fig:UMAP_and_FRC_RGG_den_025}, we show classification of geometric graphs with fixed density or $\rho=0.25$ and different box dimensions. In particular, we compute the global $d$-FRC for $d \in \{1,2,3\}$ (top panels) and parameterize UMAP with these geometric figures. The UMAP embedding is depicted in the bottom panels.  Notably, a similar embedding pattern is observed for different $d$. This is in part because the higher orders $(d=2,\,d=3)$ information builds upon pairwise interactions ($d=1$). Specifically, the distance-based cliques are computed from pairwise information. In other words, the tetrahedra are dependent on the preceding existence of triangles, which only exist from the coupling of three edges. The above results are robust against other network density values (see~\autoref{sec:RGG_figures} for further examples).
%
%and compare it with the geometric classification for %$\rho=0.25$. It is clear that the statistical differences %between dimensionality groups are enhanced by the data %geometrization performed by the FRC. Notably, the same %pattern is observed for different $d$. However, the %information from pairwise interactions ($d=1$) presents %equivalent results to higher orders $(d=2,\,d=3)$. This is %because the distance-based cliques are limited by the %evolving construction of pairwise information. Namely, the %tetrahedra are dependent on the preceding existence of %triangles, which only exist from the coupling of three edges

\begin{figure}
    \centering
    \includegraphics[width=\linewidth]{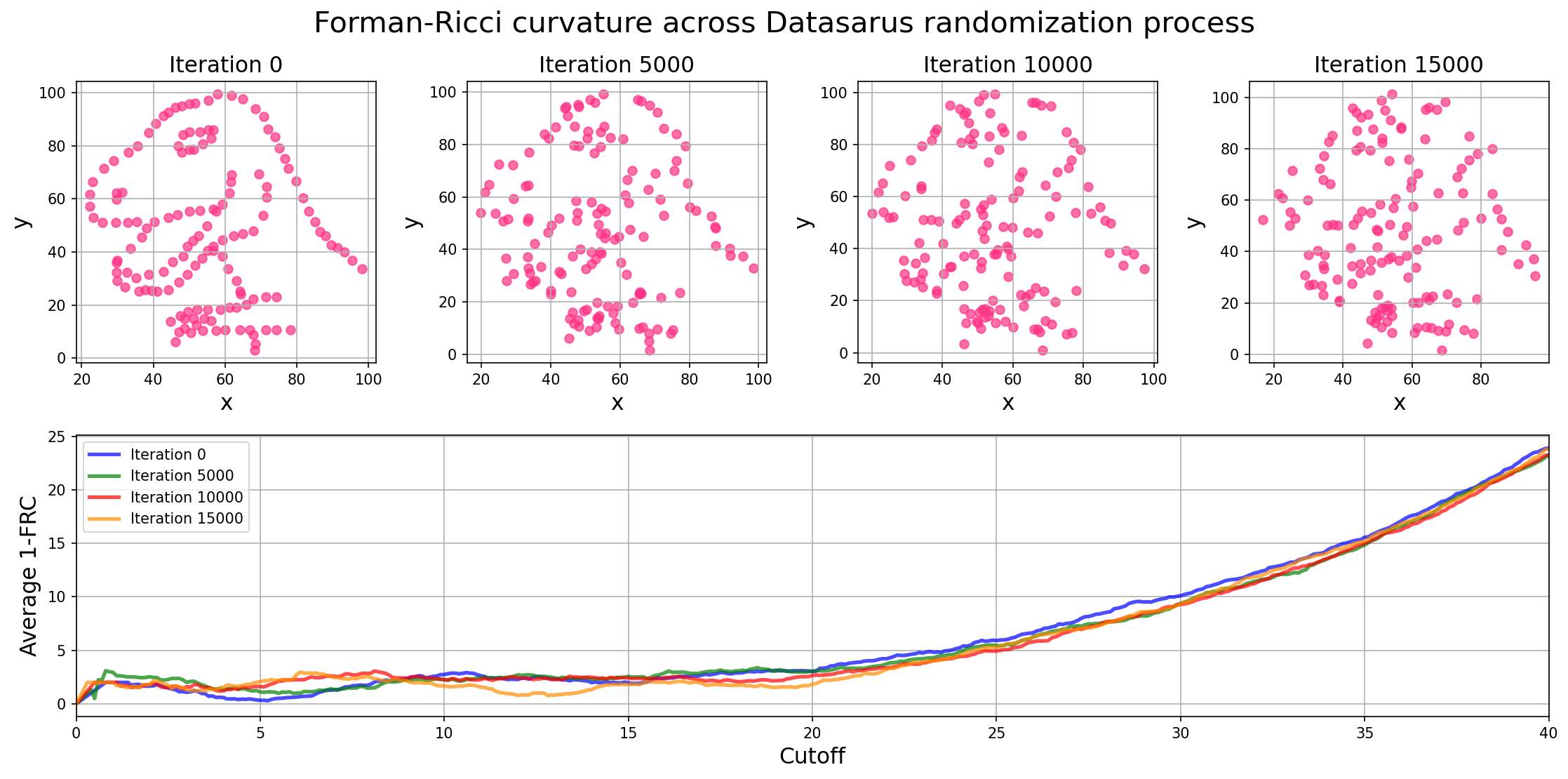}
    \caption{FRC computation for different steps of Datasaurus dataset randomization. In the example, the algorithm~\ref{alg:randomizer} was performed to disturb the original data (iteration 0) in a total of $9000$ iterations, which totalized $17469$ effective iterations. From these, we show the dataset changes for iterations $5000$, $10000$ and $15000$. Despite the randomization process that generates distinct geometry from VR complexes, the FRC presents robustness for low noise levels.}
    \label{fig:FRC_and_datasaurus_itertions}
\end{figure}

Moving to the case of DataSaurus datasets, we also test the sensitivity to noise added from the randomization iteration process. \autoref{fig:FRC_and_datasaurus_itertions} illustrates the sensitivity of the global FRC over the data randomization process.
\begin{figure}
    \centering
    \includegraphics[width=\linewidth]{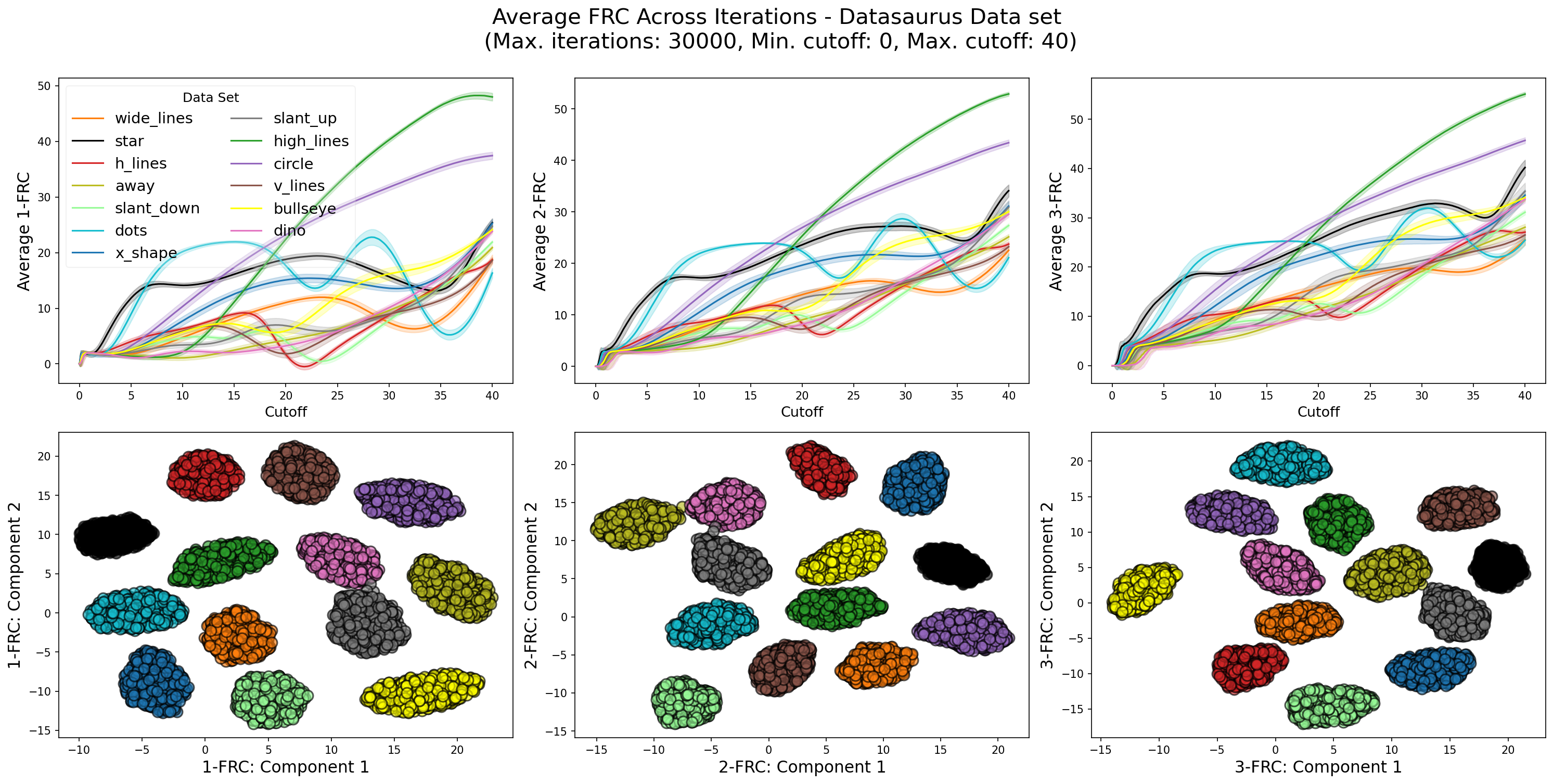}
    \caption{ Evolution of the average Forman--Ricci curvature during the Datasaurus randomization process. The first, second, and third columns show the Average $d$-dimensional Forman--Ricci curvature ($d$-FRC) for $d=1$, $d=2$, and $d=3$, respectively, computed throughout the Datasaurus randomization procedure. Each curve represents the evolution of the average $d$-FRC over a maximum of $30,000$ randomization steps. The randomized datasets were generated using Algorithm~\ref{alg:randomizer}.}
    \label{fig:FRC_and_UMAP_datasaurus_30000_iterations}
\end{figure}
In~\autoref{fig:FRC_and_UMAP_datasaurus_30000_iterations}, we provide the comparison between the geometry from Datasaurus datasets randomization and its UMAP embeddings. Noteworthy, we show that the global FRC provides a structural clustering even in the presence of noise after the randomization steps. Similar to the results on random networks, higher-order curvatures do not provide additional information since they are built upon pairwise information (i.e., due to the lack of independence of face generation in distance-based approaches). In~\autoref{sec:figures}, we test the sensitivity to higher levels of noise. Crucially, our approach provides embeddings from data geometrization of the Datasaurus dataset that generalizes the results obtained in~\cite{dlotko2023persistence}, where the authors used homology to identify statistical differences in the data.
%------------------
\subsubsection{FRC computations on breast cancer datasets}
%------------------
We now test our approach on two sets of breast cancer data, namely:

\begin{enumerate}
    \item Breast cancer diagnosis from Wisconsin dataset~\cite{UCI_Wisconsin_Breast_Cancer} with two classification labels (benign and malignant);
    \item Brest cancer classification from the Molecular Taxonomy of Breast Cancer International Consortium (METABRIC)~\cite{METABRIC_Kaggle}. We used the patient's vital status (alive, died of the disease, and died of other causes) for the classification labels. 
\end{enumerate}
In both datasets, we restricted the dataset to the numerical data features and only considered patients without missing information. We apply the geometrization process as described in \autoref{tab:Geometrized_table}, by benchmarking the geometry generated from Ollivier-Ricci curvature \cite{jost2014ollivier} (ORC  -- computed  via GraphRicciCurvature \cite{ni2024graphriccicurvature}) and Forman-Ricci curvature. We also use the silhouette score \cite{shahapure2020cluster} to benchmark the quality of the geometrical clusterings provided by each curvature For the METABRIC dataset, we additionally included a feature representing the duration of time the patient lived with the tumor, calculated as the difference between the date of death and the treatment start date. Due to computational limitations, we restricted the computation of the local FRC to edges ($1$-FRC), and also limite the benchmark between FRC and ORC to the breat cancer datasets.
\newpage
\begin{figure}[!h]
    \centering
    \includegraphics[width=\linewidth]{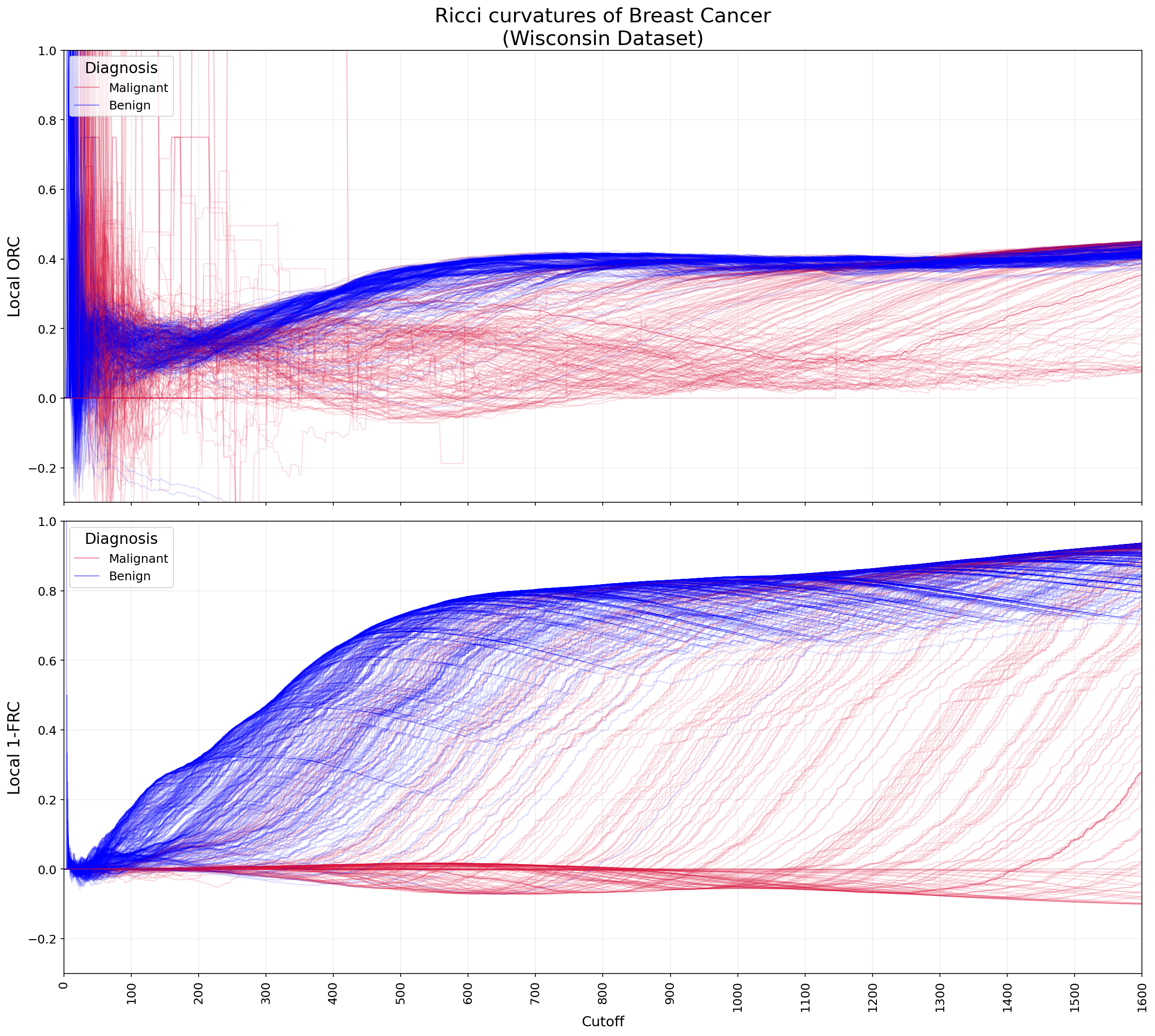}
    \caption{Local ORC (top) and 1-FRC (bottom) per patient from Wisconsin breast cancer database. The benign and malignant diagnoses are represented in blue and red, respectively. Notoriously, there is a geometric separation between the two groups from the statistical features.}
    \label{fig:FRC_bc_wisconsin}
\end{figure}
%
%\MDcom[inline]{I would remove the sentences at the top of the figures since these informations are included in the captions.}
%
\begin{figure}[!h]
    \centering
    \includegraphics[width=\linewidth]{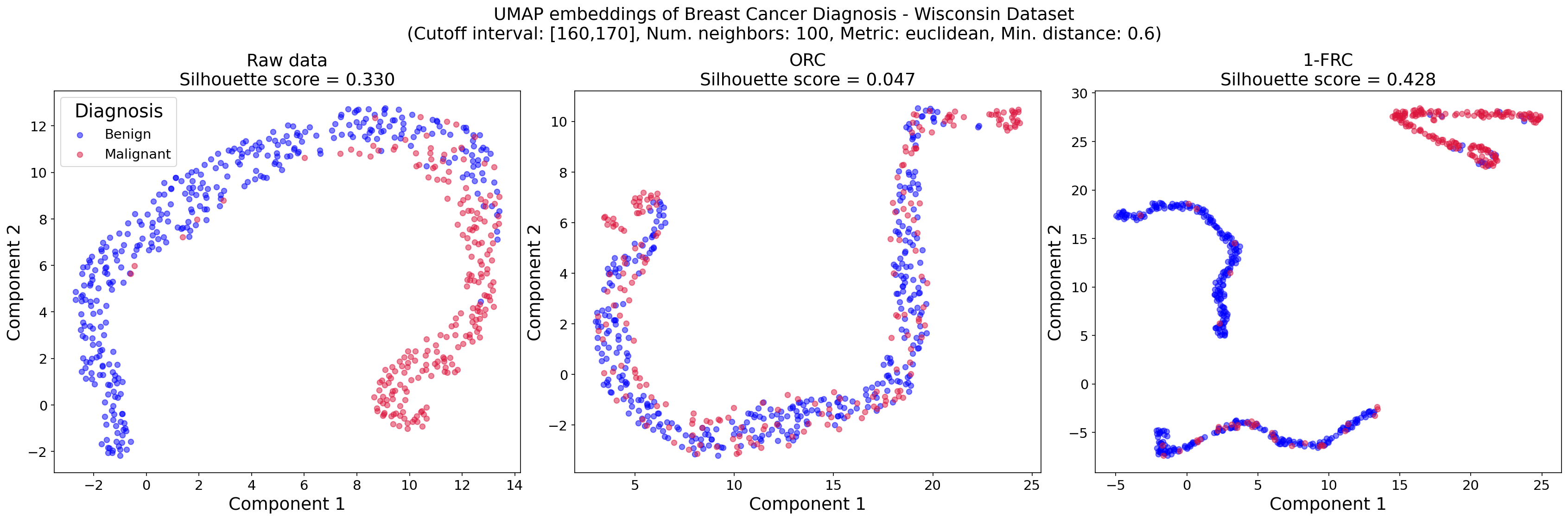}
    \caption{Comparison between the 2D-UMAP embedding of breast cancer (Wisconsin database) by using raw data (left panel)  and geometrized data {by using ORC and 1-FRC, respectively}  (\textcolor{black}{central and }right panel). {The silhouette score confirms the improvement on the clustering quality by using data geometrization.}}
    \label{fig:BC_FRC_and_UMAP_euclidean_wisconsin_160_170}
\end{figure}
\begin{figure}[!h]
%\textcolor{red}{Updated figure}
    \centering
    \includegraphics[width=\linewidth]{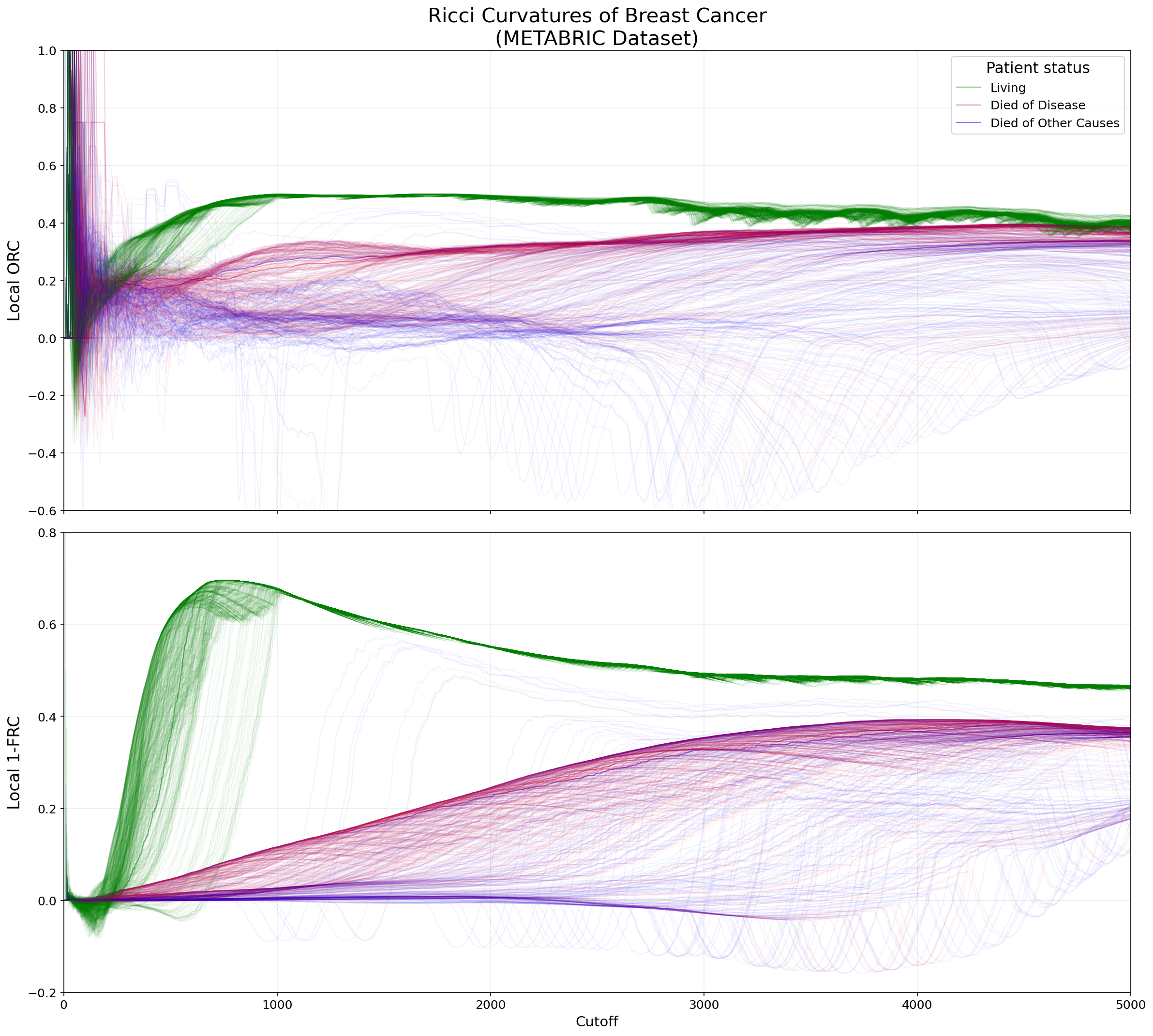}
    \caption{\textcolor{black}{Local} \textcolor{black}{ORC (top) and }1-FRC \textcolor{black}{(bottom) }per patient from METABRIC dataset. Green, blue and red curves represent the patients who survived, died of other causes and died of the disease, respectively.}
    \label{fig:FRC_bc_metabric}
\end{figure}
%
%\MDcom[inline]{I would remove the sentences at the top of the figures since these informations are included in the captions.}
%
\begin{figure}
    \centering
    \includegraphics[width=\linewidth]{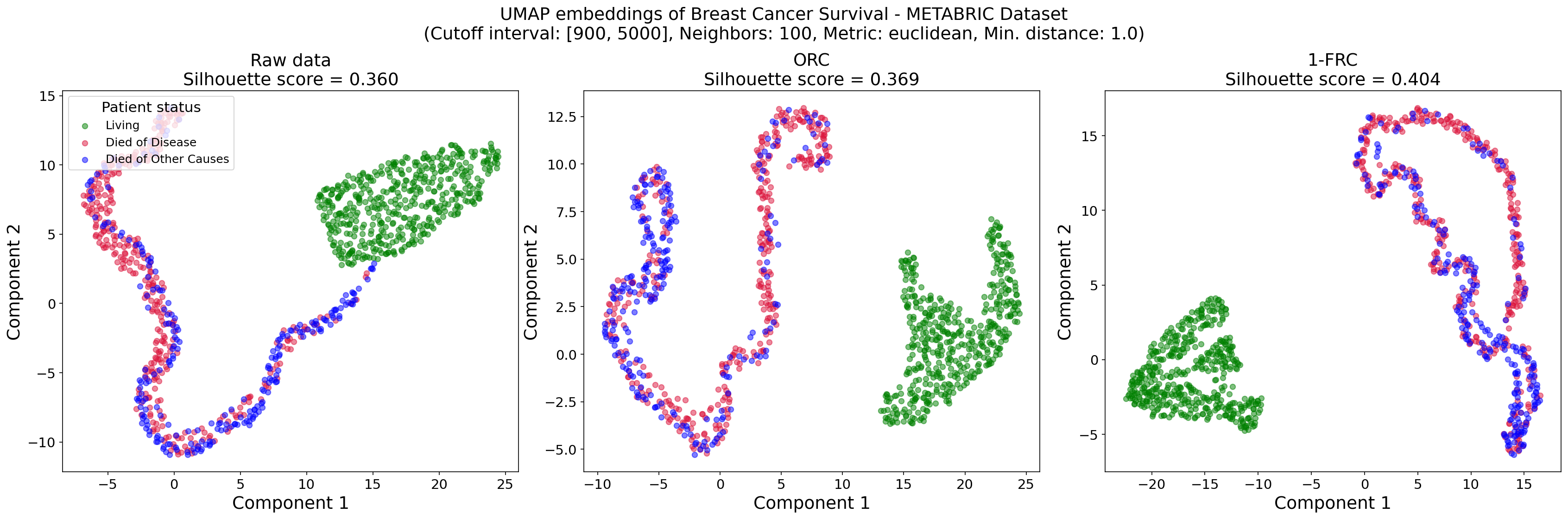}
    \caption{Comparison between the UMAP embeddings of breast cancer (METABRIC database) by using raw data (left panel)  and geometrized data by using ORC and 1-FRC, respectively (central and right panel). Notably, there is a more evident geometric separation between the groups (which is confirmed by the silhouette score improvement).}
    \label{fig:BC_FRC_and_UMAP_euclidean_METABRIC_7000_9000}
\end{figure}
In \cref{fig:FRC_bc_wisconsin,fig:FRC_bc_metabric}, we show the local $1$-FRC as a function of the cutoff distance for each patient from the Wisconsin database and METABRIC database, respectively. The resulting classification is represented via different colorings of the FRC curves. In the Wisconsin data, we used red for the malignant tumors and blue for the benign ones, respectively. For the METABRIC data, we used green, blue and red to distinguish the groups of alive, and died of other causes and died of the diseases, respectively. 

The FRC can be useful to reveal patterns embedded in the dataset's intrinsic geometry. However, it is important to emphasize that data geometrization should be viewed as a toolkit—a method that guides analysis rather than a one-time procedure. It should be approached as an iterative process, informed by continuous visualization of the resulting classifications. To illustrate this, consider a dataset composed of statistical features. Applying data geometrization to these features helps identify optimal cutoff intervals where statistical differences between data points become apparent, and it highlights statistical outliers. This insight can be directly leveraged to inform clustering decisions and improve the interpretability of the results. For example, in \autoref{fig:FRC_bc_wisconsin}, the FRC becomes less distinguishable for cutoff values above $2500$, as most curvature values converge toward similar magnitudes. This convergence suggests that the corresponding data points are not statistically separable, indicating that cutoff values greater than $2500$ are suboptimal for distinguishing between malignant and benign tumors. Conversely, cutoff intervals between $50$ and $2000$ reveal clear statistical divergence between the two groups, making them suitable for classification purposes. Similarly, in~\autoref{fig:FRC_bc_metabric}, the choice of cutoff interval plays a critical role in balancing outlier detection, redundancy reduction, and  clustering quality.

To enhance the result, we use UMAP to compare the usual data embedding (from the data features) and the FRC as a function of the cutoff interval (geometrized data). The results are shown in~\cref{fig:BC_FRC_and_UMAP_euclidean_wisconsin_160_170,fig:BC_FRC_and_UMAP_euclidean_METABRIC_7000_9000}. A key advantage of our approach is the ability to select a cutoff interval that reveals the greatest statistical divergence between groups. Although the cutoff can be chosen arbitrarily, a more systematic approach could be employed to avoid relying exclusively on visual geometric separation. For example, the silhouette score could be used to identify the optimal cutoffs that yields optimal statistical separation between groups and consequently, the optimal clustering. In \autoref{sec:figures_breast_wisconsin} and \ref{sec:breast_cancer_metabric_figures} we test the sensitivity between other cutoff intervals and different UMAP input parameters, e.g., the metric, the minimum distance between points, and the number of neighbours.
%
%-------------------------
\section{Conclusion}
\label{sec:conclusion}
%-------------------------
In this work, we introduced an alternative set-theoretical algorithm for computing local and global high dimensional Forman-Ricci curvature (FRC) in Vietoris-Rips complexes, based on the numerical increments of local curvature values as a function of cutoff distances. We also provided an alternative theoretical formulation for global geometric computation from local curvature definitions, drawing parallels with topological approaches like the Gauss-Bonnet theorem in complex networks. To benchmark our approach, we applied both global and local FRC to synthetic data for testing noise sensitivity and further compared with Olliveir-Ricci curvature applied to real-world datasets. We refer to such as method as "data geometrization". We then applied UMAP algorithm  for data dimensionality reduction, and finally benchmarked the clustering quality via silhouette score.

In contrast to state-of-the-art data dimensionality reduction methods, our data geometrization approach demonstrates robustness against noise and enhances the results of existing dimensionality reduction techniques, such as UMAP. This improvement stems from the fact that existing geometry-based dimensionality reduction algorithms assume data uniformity for accurate embeddings. Since real-world data typically lacks such uniformity, our data geometrization process offers an alternative way to find the optimal balance between signal and noise, leading to improved accuracy. This is achieved by selecting appropriate cutoff distance intervals within the Vietoris-Rips (VR) complex.

However, despite these improvements, a visual inspection of the FRC is still necessary, which can become exhaustive when dealing with large datasets containing multiple labels. Therefore, a qualitative evaluation of different cutoff values (e.g., using the silhouette score) should be performed to identify the optimal parameter space. A systematic investigation of this optimization will be the subject of future work. Additionally, the geometric computations involved remain time-consuming, limiting their application to low-dimensional structures and requiring significant computational resources and efficient code implementation. In our tests, the FRC filtration was computed using Euclidean distance (due to a limitation in the Gudhi algorithm), but the approach can be extended to other metrics, requiring alternative clique algorithms capable of constructing VR complexes in different metric spaces.

The results show that the FRC enhances the statistical relevance of clustering groups, as it allows the improvement of the sillouete score from the embeddings generated from geometrized data. These findings not only offer an alternative to barcode homology representations of datasets but also improve data embedding outcomes compared to traditional methods. In conclusion, our approach facilitates the statistical analysis of large datasets through a geometric lens.
\newpage
\bibliographystyle{siamplain}
\bibliography{references}

\newpage
section{Supplementary Materials}\label{sec:sup}

Herein, we provide the theoretical results used to derive the set-theoretical formulations that extend the computation of FRC to VR complexes. We also provide the algorithms and the proofs for our theoretical findings. We also include a few additional figures.

%-----------------------------------------
\section{Theoretical Results}
\label{sec:previous_results}
%-----------------------------------------
Below, we provide a set-theoretical representation of simplicial complexes, and we give an alternative formulation to the FRC. This derivation used the disjoint union $N_\alpha = P_\alpha \sqcup T_\alpha$ and the result $|T_\alpha|=(d+1)|H_\alpha|$, where $N_\alpha$, $T_\alpha$ and $P_\alpha$ are the sets of neighbours, transverse neighbours and parallel neighbours to $\alpha$, respectively. Using elementary algebra, we simplify the computation for FRC to
\begin{eqnarray}
\label{eq:previous_formula_FRC_demonstration}
\mathrm{F}(\alpha)=(d+2)\cdot|H_{\alpha}| + 2\cdot(d+1) - |N_\alpha|,
\end{eqnarray}
where $N_\alpha$ is the set of neighbors of $\alpha$.
We also obtain that
\begin{eqnarray}\label{eq:N_set_theoretic}
   N_\alpha=\bigsqcup_{\substack{\gamma\in\partial(\alpha) \\ x\in\pi_\gamma\neq \emptyset\\x\notin\alpha}}\{\gamma\cup\{x\}\},
\end{eqnarray}
 and hence, we can rewrite $|N_\alpha|$ as
\begin{eqnarray}\label{eq:new_neighborhood}
    |N_\alpha|=\sum_{\gamma \in \partial \alpha}|\pi_\gamma|-(d+1). 
\end{eqnarray}
Finally we derive the number of $(d+1)$-cells containing $\alpha$ as
\begin{eqnarray}\label{eq:higher_order_cells}
|H_\alpha|=|\bigcap_{\gamma \in \partial \alpha}\pi_\gamma|.
\end{eqnarray}
Equations~\eqref{eq:new_neighborhood} and~\eqref{eq:higher_order_cells} together with~\eqref{eq:previous_formula_FRC_demonstration} provide a new formulation for FRC computation, namely:
\begin{eqnarray}\label{eq:new_formula_FRC_demonstration}
    \mathrm{F}(\alpha)=(d+2)\cdot|\bigcap_{\gamma \in \partial \alpha}\pi_\gamma|+2\cdot (d+1)-\sum_{\gamma \in \partial \alpha}|\pi_\gamma|.
\end{eqnarray}
This last formulation is crucial in the implementation of our novel FRC algorithm as a function of cutoff distance.
\newpage
%------------------------
\section{Algorithms}
\label{sec:algorithms}
%------------------------
We now give details on the algorithms used in our approach. The Algorithm SM2.1 is a simplified (basic) version of the proposed algorithm for computing FRC in VR complexes. The~\cref{alg:FRC_raw} computes FRC in the function of the distance for VR complexes, and provides a more detailed version, which includes variable declarations and interpolation in the post-processing data. 

The~\cref{alg:randomizer} was used to randomize $n$-dimensional point cloud data so that the original statistics of the points are maintained (with an error of $0.1$). In particular, it was used to randomize the Datasaurus dataset (see~\cref{fig:datasaurus_dataset}), where $n=2$, $\text{scale}=0.5$ and $\text{temp}=1$. It is worth noting that the algorithm performs several steps of point randomization in its iterations, however, they are not always effective as the randomization is not performed when the statistical conditions are not reached in that specific iteration. Therefore, a high number of total steps must be performed in order to effectively randomize the data. For instance, in \autoref{fig:FRC_and_datasaurus_itertions}, a total of $90000$ iterations were performed to obtain that $17469$ sequential iterations provide effective data randomization.

\begin{algorithm}
\caption{Compute Average global Forman-Ricci Curvature}
%\label{alg:_tot_avg_frc}
\label{alg:FRC_raw}
\begin{algorithmic}
\STATE{Input: $C,d_{\max}\geq 1$}
%\STATE{$\mathcal{F}_d=\emptyset,\forall d\leq d_{\max}$}
\STATE{Set $\mathrm{F}_d(x):=0, \,\forall x \in V$}
\STATE{Set $c_d:=0,\forall d\in \{1,\hdots,d_{\max}\}$}
%\STATE{Set $\mathrm{F}_d(\alpha):=0, \,\forall \alpha \in C_d$}
\STATE{Set $w=\infty$}
\STATE
%\STATE{Compute \text{Neigh(G)} from \cref{alg:Neig};}
%\FOR{$d\leq d_{\max}$}
\FOR{$\alpha\in C$}
\STATE{Set $d:=|\alpha|-1$}
\STATE{$c_d:=c_d+1$}
\STATE{Update the neighborhood of each $\gamma \in \partial(\alpha)$}
\STATE{Compute $w_\alpha$ the diameter of $\alpha$}
\STATE{Compute $\mathrm{F}_{d}(\alpha)$ and $\mathrm{F}_d(x)$ according to \eqref{eq:new_formula_FRC} and the help of \eqref{eq:delta_function_node}, 
\textit{i.e.,}:
\STATE{Set the number of neighbors of $\alpha$, n:=0}
    %\item set $p:=0,\,t:=0;$
    %\item compute $\partial(\alpha)$ and $\pi_\alpha$;
    %\STATE{Compute $\pi_{\alpha};$}
    \STATE{Compute $H:=\bigcap_{\substack{\gamma \in \partial(\alpha)}} \pi_\gamma$}
    \FOR{ $\gamma \in \partial(\alpha)$}
    \STATE{update $n:=n+|\pi_\gamma|$}
    \ENDFOR
   \STATE 
    \STATE{Set $\mathrm{f}_d:=\mathrm{F}_{d}(\alpha)=(d+2)|H|+2\cdot (d+1)-n$}

}
\STATE
\FOR{$x \in V$}
\STATE{Update the contribution of $\mathrm{F}(\alpha)$ to the node $x$, \textit{i.e., }$\mathrm{F}_d(x):=\mathrm{F}_d(x)+\mathrm{f}_d/(d+1)$}
\ENDFOR
\FOR{$\gamma \in \partial(\alpha)$}
\FOR{$x \in \pi_\gamma\setminus \alpha$}
\STATE{$\delta:=\delta(x)$}
\FOR{$y \in \partial(\alpha)$}
\STATE{Update local total FRC for the nodes in the boundary: $\mathrm{F}_d(y):=\mathrm{F}_d(y)+\delta/(d+1)$}
\ENDFOR
\STATE{Update the total FRC: $\mathrm{f}_d:=\mathrm{f}_d+\delta$}
\ENDFOR
\ENDFOR
\IF{$w\neq w_\alpha$}
\STATE{}
\IF{$c_d\neq 0$}
\STATE{\textbf{print} $(w,\frac{f_{d}}{c_d},(\frac{\mathrm{F}_d(x)}{c_d},\, \forall x \in V))$}
\ELSE
\STATE{\textbf{print} $(w,0,(0,\,\forall x \in V))$}
\ENDIF
\ENDIF
%\STATE{\textbf{If} $w\neq w_\alpha$, then \textbf{print} $(w,\frac{F_{d}(\alpha)}{c_d},\frac{\mathrm{F}_d(x)}{c_d},\,\forall x \in V)$}
\STATE{Update $w:=w_\alpha$}
%\STATE{Update $\mathcal{F}_d := \mathcal{F}_d\cup \{ \mathrm{F}_{d}(\alpha)\}$}
%\ENDFOR
\ENDFOR
%\RETURN $\mathcal{F}=\bigcup\limits_{d}\{\mathcal{F}_d\}$
\end{algorithmic}
\end{algorithm}

\begin{algorithm}
\caption{Dataset Perturbation}
\label{alg:randomizer}
\begin{algorithmic}[1]
\STATE \textbf{Input:} Initial dataset $\text{initial\_ds}$, Number of iterations $\text{iterations}$, Temperature $\text{temp}$ , Perturbation $\text{scale}$.
\STATE \textbf{Output:} Final dataset $\text{current\_ds}$ with preserved statistical properties

\STATE

\STATE Set $\text{current\_ds} \leftarrow \text{initial\_ds}$

\STATE

\STATE \textbf{Function} \textsc{FIT}$(\text{ds})$:
\STATE \hspace{1em} \textbf{Return} the sum of distances from the origin for all points in $\text{ds}$

\STATE

\STATE \textbf{Function} \textsc{ISERROROK}$(\text{test\_ds}, \text{initial\_ds})$:
\STATE \hspace{1em} Calculate the mean and standard deviation for each coordinate in $\text{initial\_ds}$ and $\text{test\_ds}$
\STATE \hspace{1em} Calculate the correlation matrices for $\text{initial\_ds}$ and $\text{test\_ds}$
\STATE \hspace{1em} Round all values to 3 decimal places
\STATE \hspace{1em} \textbf{Return} True if means, standard deviations, and correlations match to 3 decimal places; otherwise, return False

\STATE

\STATE \textbf{Function} \textsc{MOVERANDOMPOINTS}$(\text{ds})$:
\STATE \hspace{1em} Create a copy of $\text{ds}$ as $\text{test\_ds}$
\STATE \hspace{1em} Select a random index $\text{idx}$ from $\text{ds}$
\STATE \hspace{1em} Generate a random movement vector from a normal distribution with small-scale
\STATE \hspace{1em} Move the point at $\text{test\_ds}[\text{idx}]$ by the movement vector
\STATE \hspace{1em} \textbf{Return} $\text{test\_ds}$

\STATE

\STATE \textbf{Function} \textsc{RANDOM}(scale):
\STATE \hspace{1em} \textbf{Return} a random value between 0 and 1 (with standard deviation = scale)

\STATE

\STATE \textbf{Function} \textsc{PERTURB}$(\text{ds}, \text{temp})$:
\STATE \hspace{1em} \textbf{Loop} until a valid perturbation is found:
\STATE \hspace{2em} $\text{test} \leftarrow \textsc{MOVERANDOMPOINTS}(\text{ds})$
\STATE \hspace{2em} \textbf{If} $\textsc{FIT}(\text{test}) > \textsc{FIT}(\text{ds})$ \textbf{or} $\text{temp} > \textsc{RANDOM()}$:
\STATE \hspace{3em} \textbf{Return} $\text{test}$

\STATE

\STATE \textbf{Main Function} \textsc{SimulatedAnnealing}$(\text{initial\_ds}, \text{iterations}, \text{temp})$:
\STATE \hspace{1em} Set $\text{current\_ds} \leftarrow \text{initial\_ds}$

\STATE \hspace{1em} \textbf{For} each iteration in $1$ to $\text{iterations}$:
\STATE \hspace{2em} $\text{test\_ds} \leftarrow \textsc{PERTURB}(\text{current\_ds}, \text{temp})$
\STATE \hspace{2em} \textbf{If} \textsc{ISERROROK}$(\text{test\_ds}, \text{initial\_ds})$:
\STATE \hspace{3em} Update $\text{current\_ds} \leftarrow \text{test\_ds}$

\STATE \hspace{1em} \textbf{Return} $\text{current\_ds}$

\end{algorithmic}
\end{algorithm}
\newpage

%-------------------------------------------------------------
\section{Main results and Demonstrations}
\label{sec:demonstrations}
%-------------------------------------------------------------
In this section, we develop the theoretical assumptions of our work.
\begin{prop}\label{prop:before_after}
    Let $\alpha_i \in C_d$, $\alpha \in N^i$. Let $\mathrm{F}:C_d\rightarrow \mathbb{Z}$ and $\mathrm{F}^i$ be the FRC function as defined in \ref{eq:original_FRC}. If $\alpha$ is parallel to $\alpha_i$, then $\mathrm{F}^i(\alpha)=\mathrm{F}^{i-1}(\alpha)-1$. Otherwise, we have  $\mathrm{F}^i(\alpha)=\mathrm{F}^{i-1}(\alpha)+(d+1)$.
\end{prop}
\begin{proof}
   We have that  $F^{i-1}(\alpha)=|H^{i-1}_{\alpha}|+(d+1)-|P^{i-1}_\alpha|$. Suppose that $\alpha$ is parallel to $\alpha_i$. Thus, $P^i_\alpha=P^{i-1}_{\alpha}\sqcup\{\alpha_i\}$, which implies that $|P^{i}_{\alpha}|=|P^{i-1}_{\alpha}|+1$ and that $H^{i}_\alpha=H^{i-1}_{\alpha}$. It follows that $\mathrm{F^i}(\alpha)=|H^{i}_{\alpha}|+(d+1)-|P^{i}_\alpha|=|H^{i-1}_{\alpha}|+(d+1)-|P^{i-1}_\alpha|-1=\mathrm{F^{i-1}}(\alpha)-1.$ Suppose that $\alpha \notin P^{i+1}.$ Then, $\alpha \in T^{i+1},$ which implies that $|H^{i}_\alpha|=|H^{i-1}_\alpha|+1$ and \textcolor{black}{$|P^i_\alpha|=|P^{i-1}_\alpha|-d$}. It follows that $\mathrm{F}^i_{\alpha}=|H^{i}_\alpha|+(d+1)-|P^{i}_\alpha|=$ $\textcolor{black}{(|H^{i-1}_\alpha|+1)+(d+1)-(|P^{i-1}_\alpha|-d)}=\mathrm{F}^{i-1}(\alpha)+(d+1).$
\end{proof}
\begin{theorem}[Geometric Gauss-Bonnet Theorem]\label{theo:Geometric_Gauss_Bonnet}
    Let $C_d$ be a non-empty set of $d$-faces in a simplicial complex $C$ generated from an undirected graph $G=(V,E)$. Let $F_d$ and $f_d$ be the global and local FRC definitions, as defined in \eqref{eq:avg_FRC} and \eqref{eq:local_total_forman}, respectively. Then, the following equality holds:
  \begin{eqnarray}
    \mathrm{F}_d(C)=\sum_{x\in V} \mathrm{f}_d(x).
\end{eqnarray}
\end{theorem}
\begin{proof}
    It is sufficient to notice that
    \begin{eqnarray}
        \sum_{x \in V}\sum_{\substack{\alpha \subset C_d \\ x \in \alpha}} \mathrm{F}_d(
\alpha)=(d+1)\sum_{\alpha \subset C_d}\mathrm{f}_d(\alpha),
\end{eqnarray}
once that for each $\alpha \in C_d$, the value $\mathrm{f}_d(\alpha)$ is counted exactly $(d+1)$ times in the sum over all the nodes $x \in V$. The result follows by multiplying both sides of the equation above by $\frac{1}{(d+1)\cdot |C_d|}$. %\autoref{fig:datasaurus_dataset} shows the 2D plots of the Datasaurus dataset.
\end{proof}
%
%\MDcom[inline]{The last sentence of the proof does not belong here.}
\newpage
%

%----------------------------------------
\section{Additional figures}
\label{sec:figures}
\begin{figure}[!h]
    \centering
    \includegraphics[width=\linewidth]{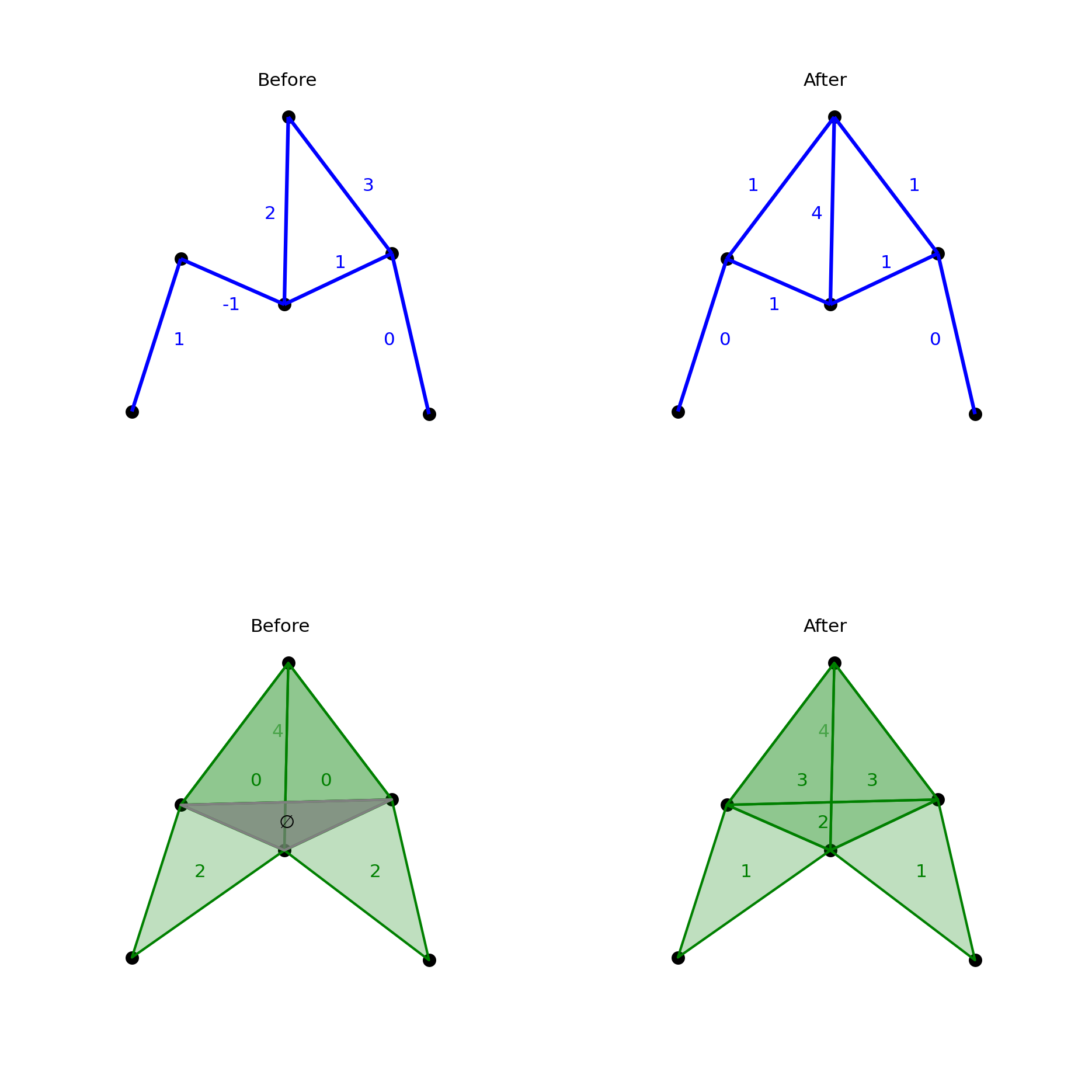}
    \caption{Example of an evolving simplicial complex (by adding a new face) and how it influences the FRC (the numerical values given in the figure) for edges (blue) and triangles (green). The new edge has $4$ neighbours, of which two are transverse (the neighbours sharing the same triangle) and $2$ parallel (The ones outside the triangle). The $1$-FRC of the new edge is the number of triangles containing the new edge plus the length of the boundary minus the number of parallel neighbours, \textit{i.e.,} $1+2-2=1.$ When a new edge is added, the FRC increases by 2 units when a new triangle is created and decreases by 1 unit for the new neighbours outside the new triangle. Similarly, the new triangle has $3$ neighbours, where $1$ shares a tetrahedron and $2$ are parallel neighbours, therefore, the $2$-FRC of the new triangle is $1+3-2=2$. Similarly, for triangle faces, the FRC increases by 3 units when a new tetrahedron is created in the neighbourhood and decreases by 1 unit otherwise.} 
    \label{fig:example_adding_face}
\end{figure}

%----------------------------------------
We now provide all the supplementary figures of our work, which include examples, illustrations and results. \autoref{fig:example_adding_face} provides an example of FRC value change as a function of the local neighbourhood. \autoref{fig:datasaurus_dataset} shows the 2D plots of the Datasaurus dataset.
%
%-------------
\subsection{Random geometric graphs}
\label{sec:RGG_figures}
\begin{figure}
    \centering
    \includegraphics[width=\linewidth]{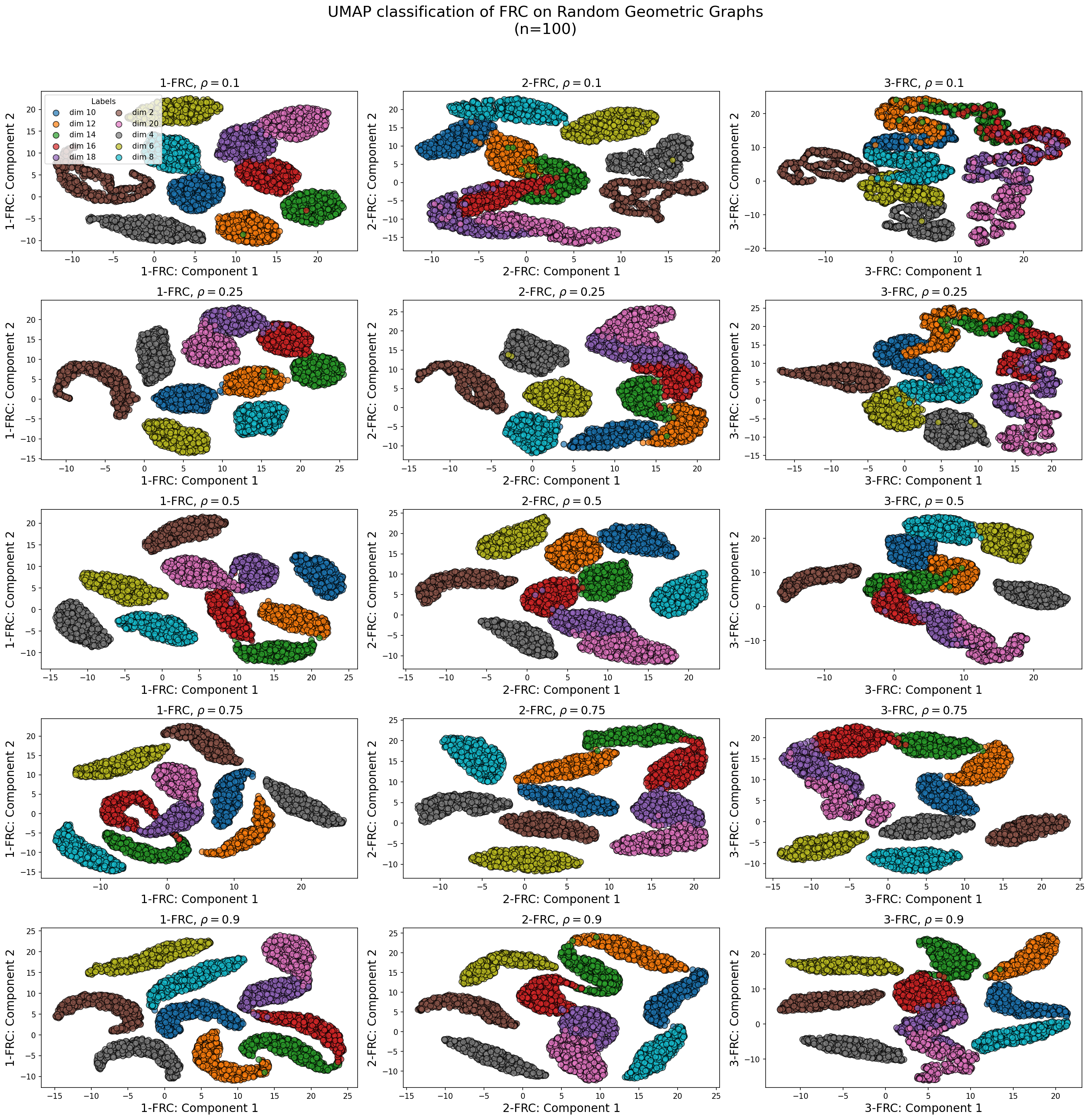}
    \caption{Computations of the average $d$-FRC for $d=1$ (left column), $d=2$ (middle column) and $d=3$ (right column) on Random geometric graphs with $n=100$ nodes and different box dimensions and densities (plot lines). The solid-coloured lines are the average, while the error bands are computed from the standard deviation. The box dimension classification can be better visualised in~\autoref{fig:UMAP_FRC_RGG_all_plots}.}
    \label{fig:FRC_RGG_all_plots}
\end{figure}

\begin{figure}
    \centering
    \includegraphics[width=\linewidth]{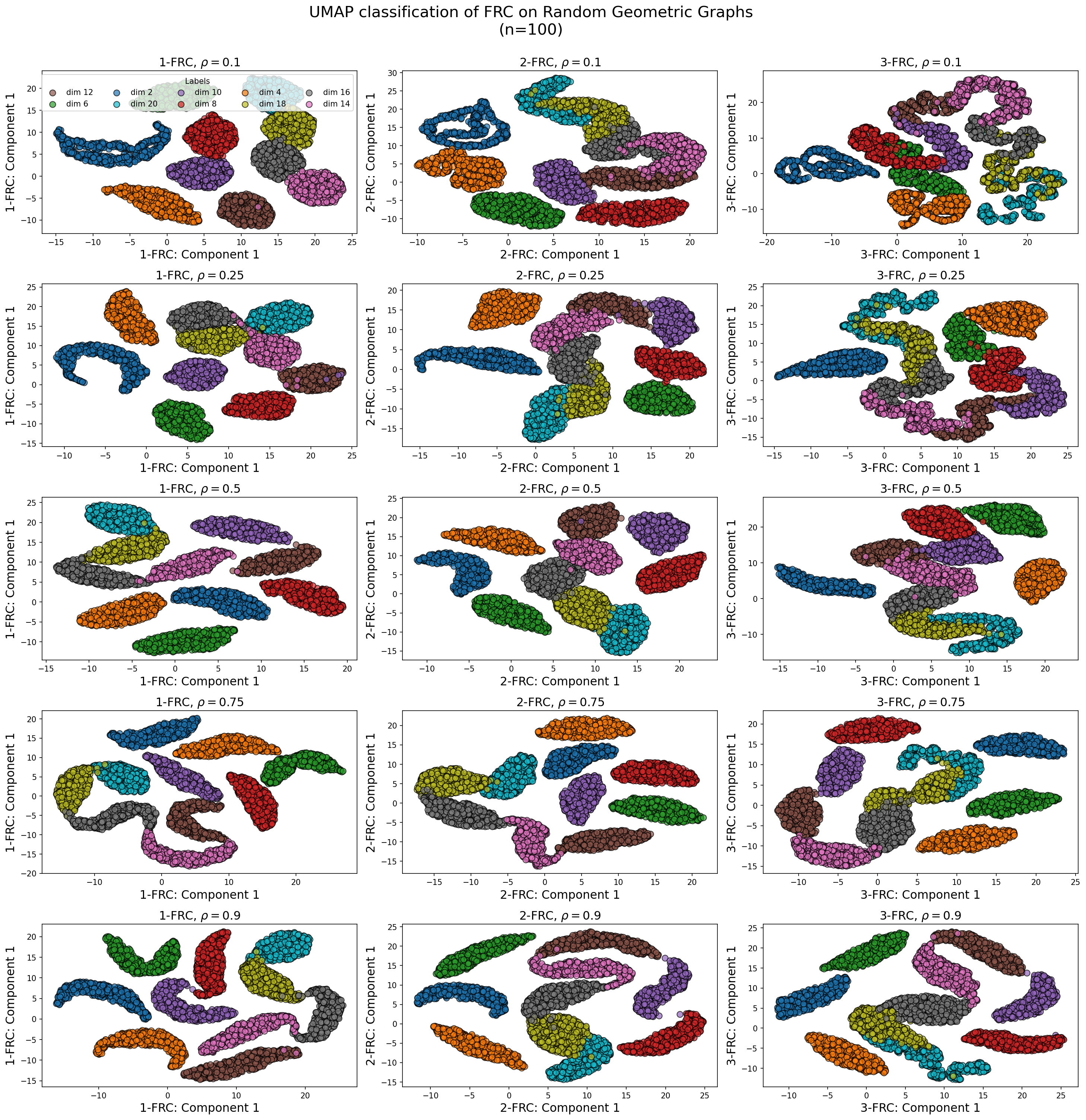}
    \caption{UMAP classification from the computations of the average $d$-FRC for $d=1$ (left column), $d=2$ (middle column) and $d=3$ (right column) on random geometric graphs with $n=100$ nodes and different box dimensions and densities; see also~\autoref{fig:FRC_RGG_all_plots}}.
    \label{fig:UMAP_FRC_RGG_all_plots}
\end{figure}

%-------------
\cref{fig:FRC_RGG_all_plots,fig:UMAP_FRC_RGG_all_plots} provide the FRC computations on random geometric graphs for different edge densities and box dimensions, as well as the comparison with the geometric classification performed by UMAP. In all tests, we used the Euclidean distance and minimum distance of $1.0$ as UMAP parameters.
%
%-------------
\subsection{Datasaurus dataset}
\label{sec:datasaurus_figures}
\begin{figure}[!h]
    \centering
    \includegraphics[width=1.\linewidth]{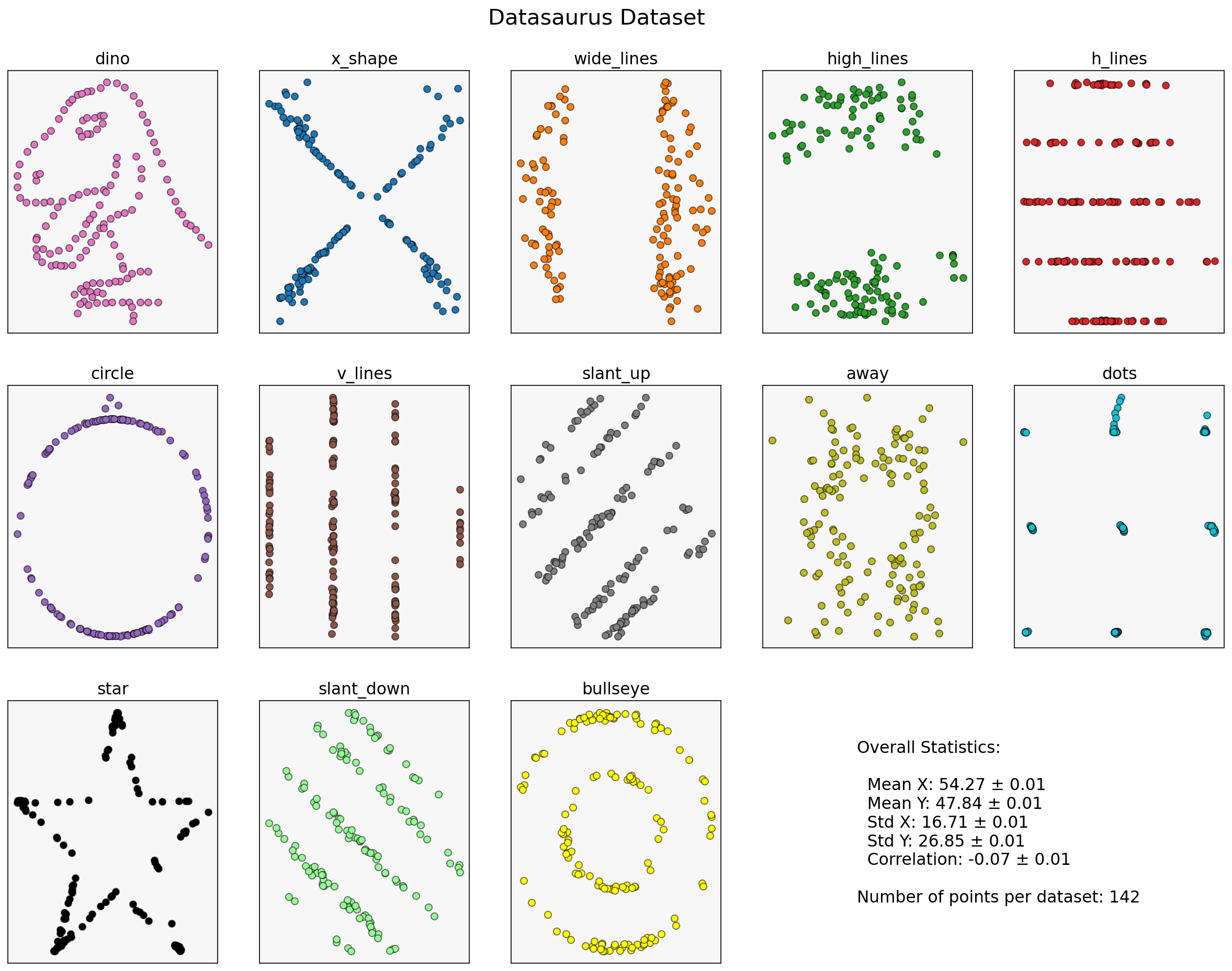}
    \caption{The Datasaurus dataset and its overall statistics. The purpose of this dataset was to generate 2D point cloud data that provides different geometries with (proximally) the same basic statistics.}
    \label{fig:datasaurus_dataset}
\end{figure}

\begin{figure}
    \centering
    \includegraphics[width=1.\linewidth]{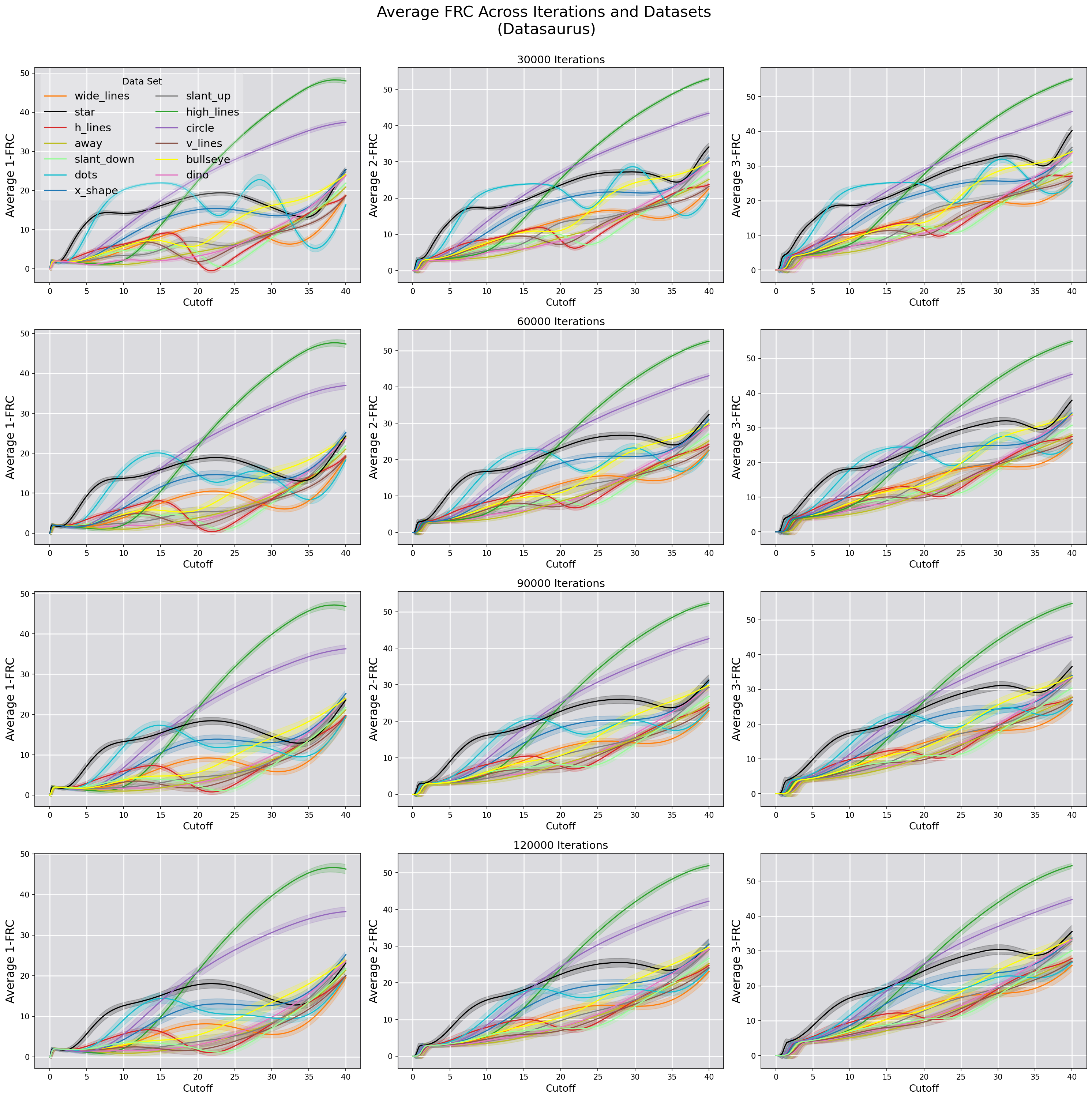}
    \caption{Average $k$-FRC statistical computation on Datasaurus datasets, for $d=1$ (left column), $d=2$ (middle column) and $d=3$ (right column), and different number of maximum randomization steps (lines). The result is the mean curvatures (central coloured lines) and the error bands (computed from the standard deviation). This result can be compared with the UMAP classification in \autoref{fig:UMAP_datasaurus_and_global_FRC_all_class}.}
    \label{fig:datasaurus_and_global_FRC_all_plots}
\end{figure}

\begin{figure}
    \centering
    \includegraphics[width=1.\linewidth]{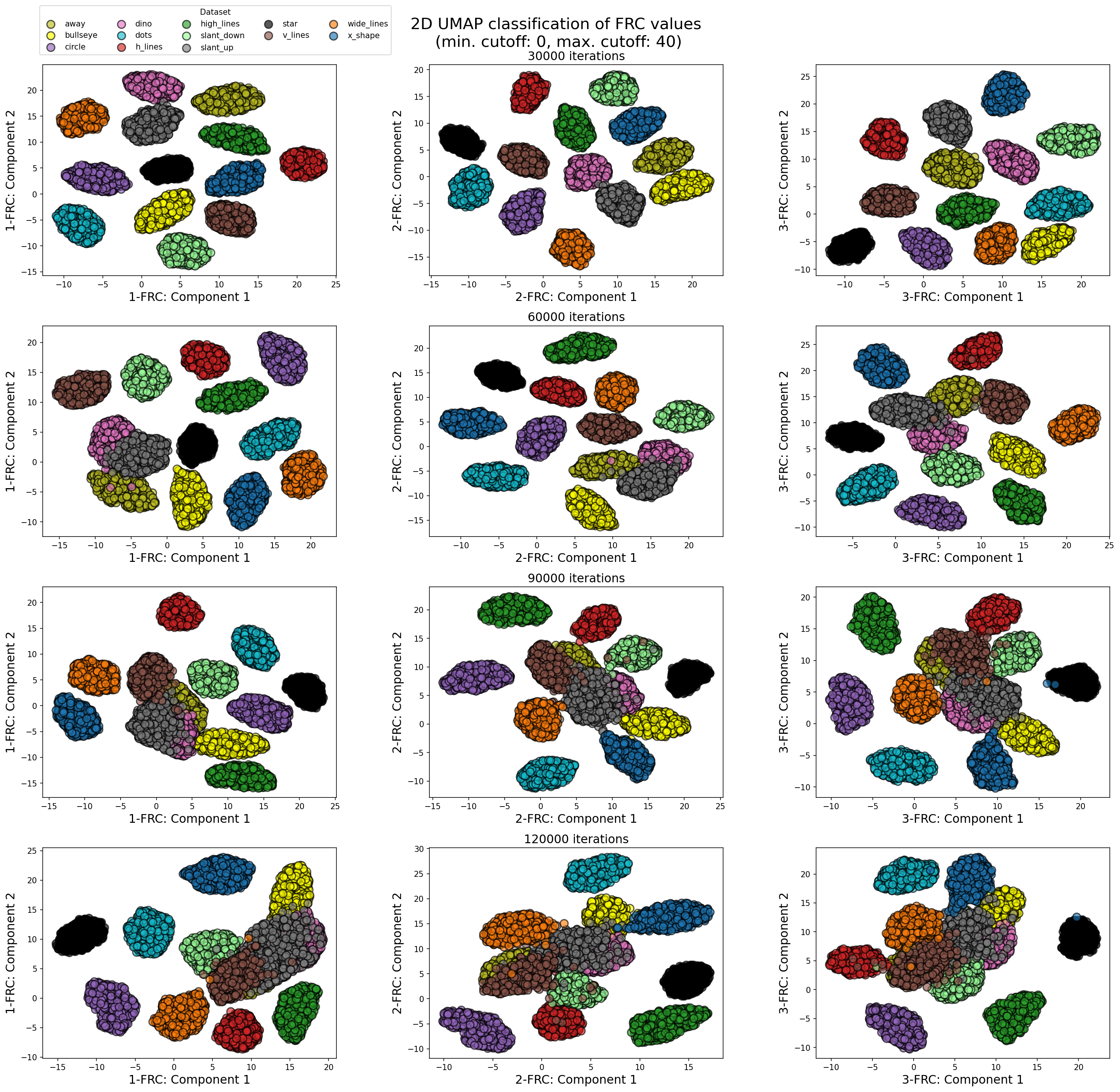}
    \caption{UMAP classification of the average $k$-FRC statistical computation on Datasaurus datasets, for $d=1$ (left column), $d=2$ (middle column) and $d=3$ (right column), and different number of maximum randomization steps (lines). Compare with}~\autoref{fig:datasaurus_and_global_FRC_all_plots}.
    \label{fig:UMAP_datasaurus_and_global_FRC_all_class}
\end{figure}

%-------------
Here, we provide the result of the data randomization process of the Datasaurus datasets performed by~\cref{alg:randomizer}, as well as its FRC computation for data classification with UMAP and sensitivity to noise. The original dataset is shown in \cref{fig:datasaurus_dataset}. In \cref{fig:datasaurus_and_global_FRC_all_plots,fig:UMAP_datasaurus_and_global_FRC_all_class} we provide the comparison of the FRC performance and classification in the presence of different levels of noise.
\newpage
%-------------
\subsection{Breast cancer Wisconsin database}
\label{sec:figures_breast_wisconsin}
\begin{figure}[!h]
%\textcolor{red}{Figure updated}
    \centering
    \includegraphics[width=1.\linewidth]{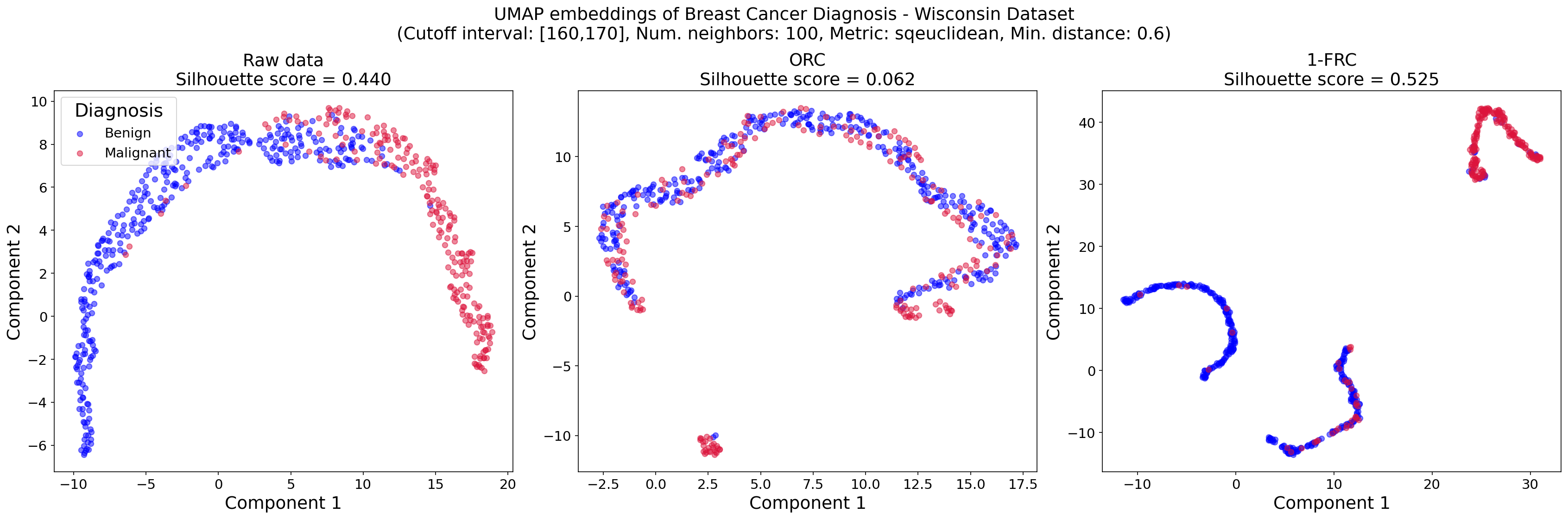}
    \caption{Breast cancer diagnosis comparison between UMAP embeddings and raw data input from Wisconsin dataset (left) and its geometrized version\textcolor{red}{s }(\textcolor{black}{center and } right).}
    \label{fig:bc_n100_md_060_c_160_170_sqeuclidean}
\end{figure}

\begin{figure}[!h]
%textcolor{red}{Figure updated}
    \centering
    \includegraphics[width=\linewidth]{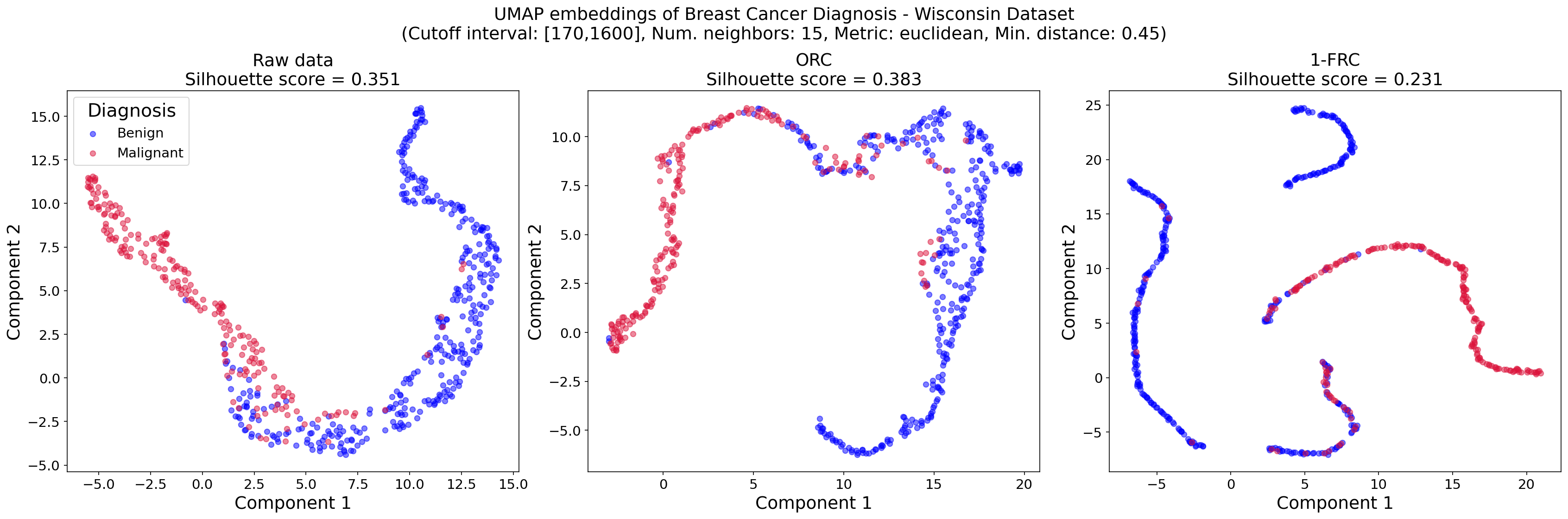}
    \caption{Breast cancer diagnosis comparison between UMAP embeddings and raw data input from Wisconsin dataset (left) and its geometrized versions(\textcolor{black}{center and }right).}
    \label{fig:bc_n15_md_045_c_170_3650_chebyshev}
\end{figure}

\begin{figure}[!h]
%\textcolor{red}{Figure updated}
    \centering
    \includegraphics[width=\linewidth]{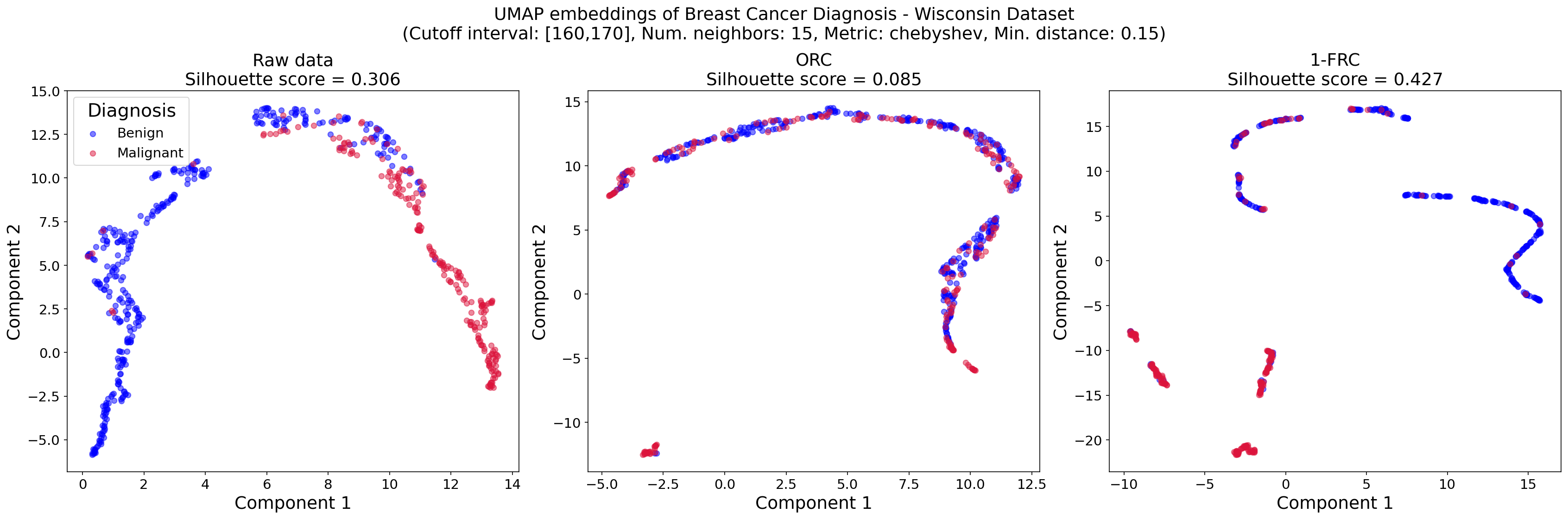}
    \caption{Breast cancer diagnosis comparison between UMAP embeddings and raw data input from Wisconsin dataset (left) and its geometrized versions (\textcolor{black}{center and} right).}
    \label{fig:bc_n15_md_015_c_160_170_chebyshev}
\end{figure}

\begin{figure}[!h]
%\textcolor{red}{Figure updated}
    \centering
    \includegraphics[width=\linewidth]{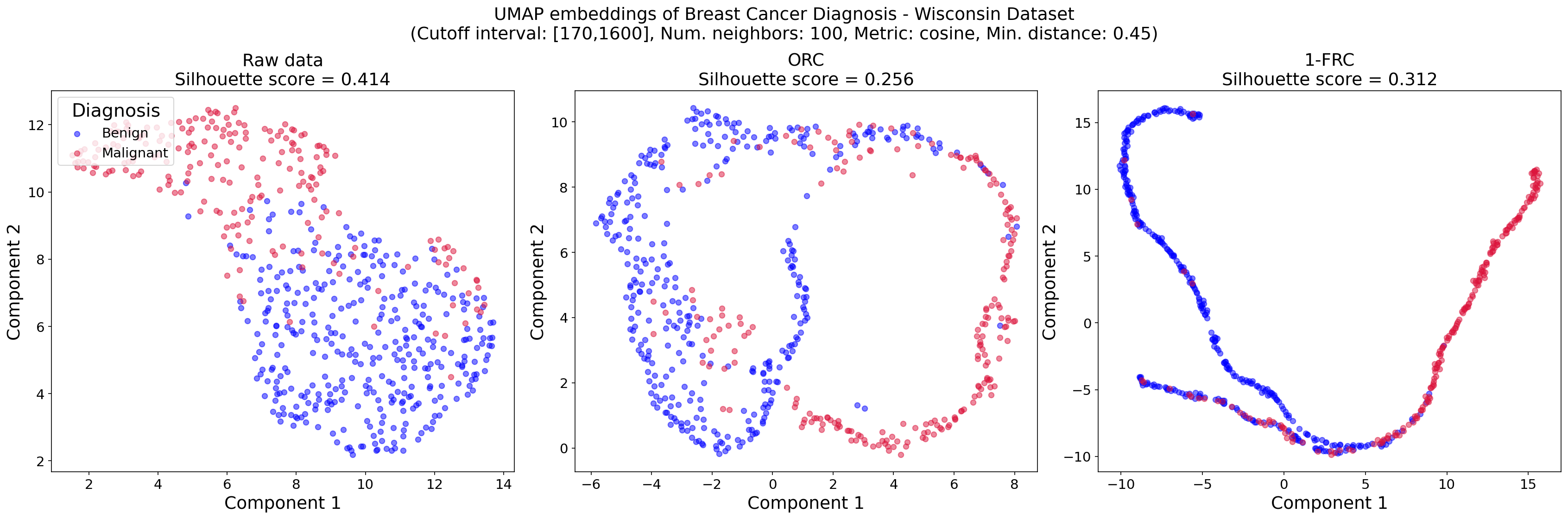}
    \caption{Breast cancer diagnosis comparison between UMAP embeddings and raw data input from Wisconsin dataset (left) and its geometrized versions (\textcolor{black}{center and }right).}
    \label{fig:bc_n100_md_045_c_170_3650_cosine}
\end{figure}

%-------------
Here, we provide the FRC computations on the Wisconsin Breast Cancer database. \cref{fig:bc_n100_md_060_c_160_170_sqeuclidean,fig:bc_n15_md_045_c_170_3650_chebyshev,fig:bc_n15_md_015_c_160_170_chebyshev,fig:bc_n15_md_015_c_160_170_chebyshev,fig:bc_n100_md_045_c_170_3650_cosine} show FRC computations and their sensitivity to cutoff intervals and UMAP parameters.
\newpage
%-------------
\subsection{Breast cancer METABRIC dataset}
\label{sec:breast_cancer_metabric_figures}
\begin{figure}[!h]
%\textcolor{red}{Figure updated}
    \centering
    \includegraphics[width=\linewidth]{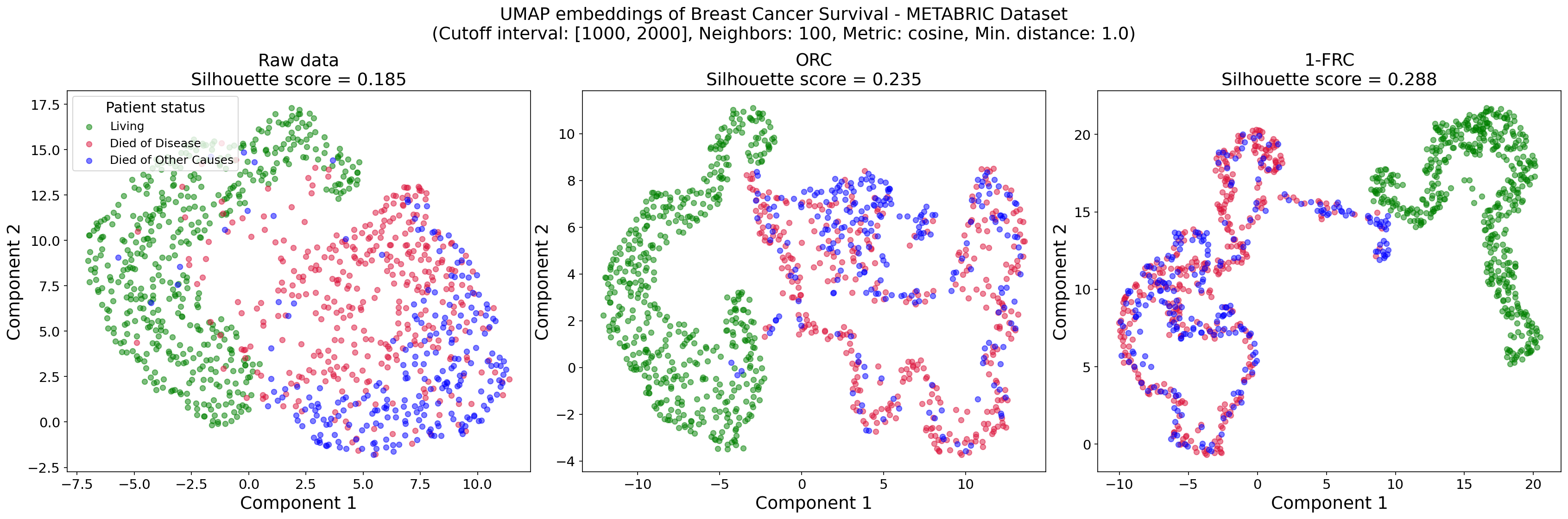}
    \caption{Breast cancer diagnosis comparison between UMAP embeddings and raw data input from METABRIC dataset (left) and its geometrized version (\textcolor{black}{center and }right).}
    \label{fig:bc_metabric_n100_md_1_c_7000_9000_cosine}
\end{figure}

\begin{figure}[!h]
%\textcolor{red}{Figure updated}
    \centering
    \includegraphics[width=\linewidth]{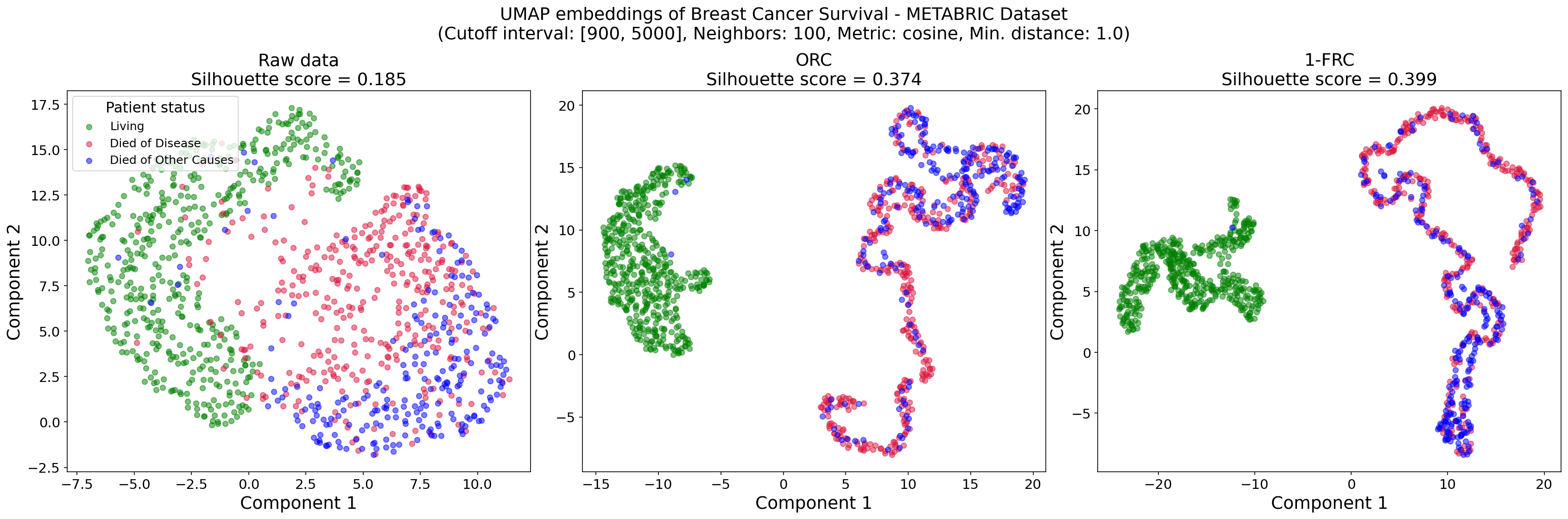}
    \caption{Breast cancer diagnosis comparison between UMAP embeddings and raw data input from METABRIC dataset (left) and its geometrized versions (\textcolor{black}{center and }right).}
    \label{fig:bc_metabric_n100_md_1_c_900_5000_cosine}
\end{figure}

\begin{figure}[!h]
%\textcolor{red}{Figure updated}
    \centering
    \includegraphics[width=\linewidth]{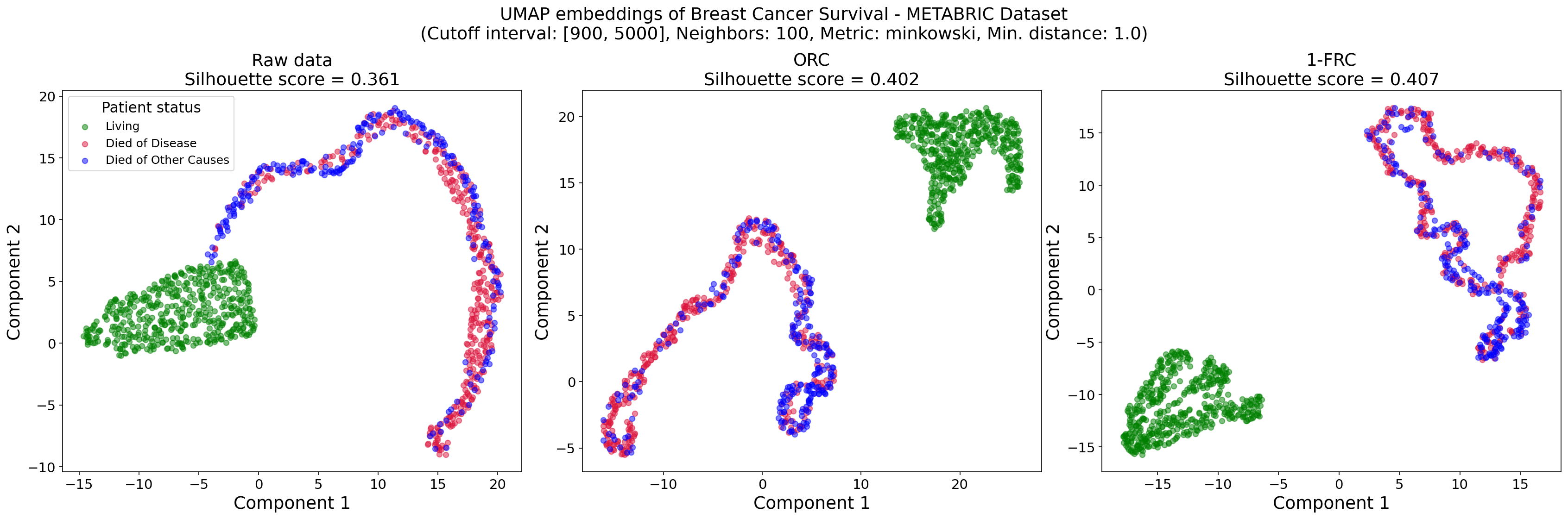}
    \caption{Breast cancer diagnosis comparison between UMAP embeddings and raw data input from METABRIC dataset (left) and its geometrized versions (\textcolor{black}{center and }right).}
    \label{fig:bc_metabric_n100_md_1_c_7000_9000_minkowski}
\end{figure}

\begin{figure}[!h]
%\textcolor{red}{Figure updated}
    \centering
    \includegraphics[width=\linewidth]{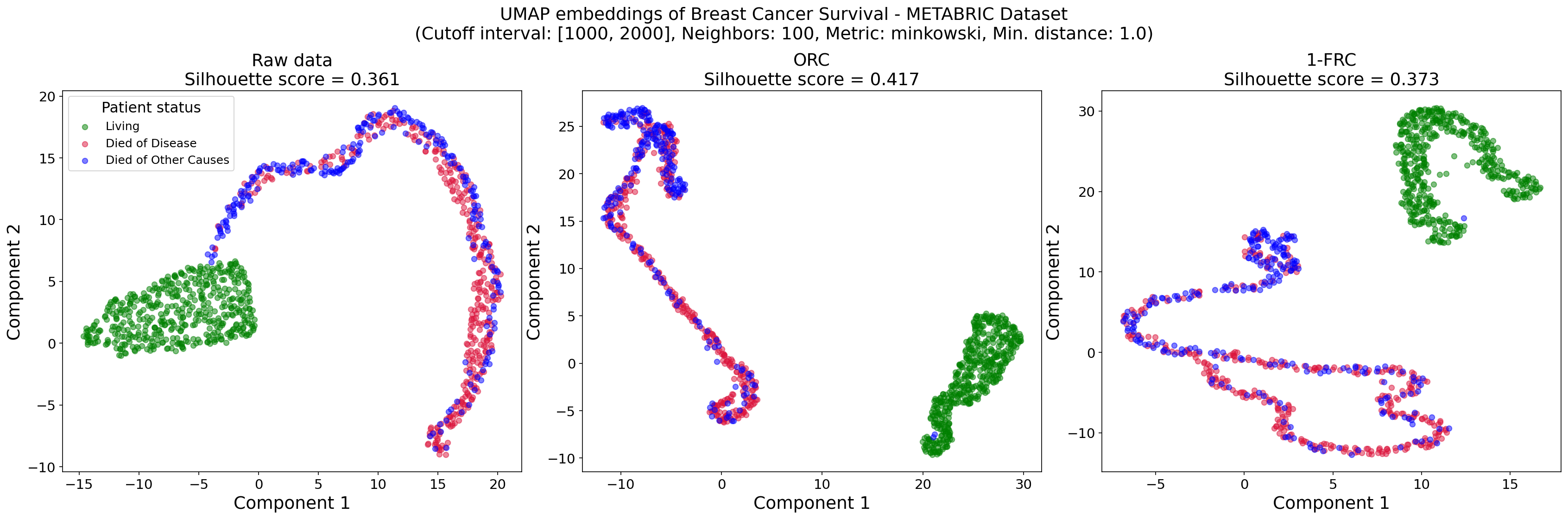}
    \caption{Breast cancer diagnosis comparison between UMAP embeddings and raw data input from METABRIC dataset (left) and its geometrized versions (\textcolor{black}{center and} right).}
    \label{fig:bc_metabric_n100_md_1_c_900_5000_minkowski}
\end{figure}

%-------------
Here, we provide the FRC computations on the breast cancer METABRIC database. \cref{fig:bc_metabric_n100_md_1_c_7000_9000_cosine,fig:bc_metabric_n100_md_1_c_900_5000_cosine,fig:bc_metabric_n100_md_1_c_7000_9000_minkowski,fig:bc_metabric_n100_md_1_c_900_5000_minkowski} display FRC computations and its sensitivity to cutoff intervals and UMAP parameters.

\end{document}

%% file: ARXIV_ex_shared.tex
\usepackage{lipsum}
\usepackage{amsfonts}
\usepackage{graphicx}
\usepackage{epstopdf}
\usepackage{algorithmic}
\ifpdf
  \DeclareGraphicsExtensions{.eps,.pdf,.png,.jpg}
\else
  \DeclareGraphicsExtensions{.eps}
\fi

\usepackage{xcolor,hyperref}
\hypersetup{colorlinks,breaklinks,
            linkcolor=blue,urlcolor=blue,
            anchorcolor=blue,citecolor=blue}

\usepackage{todonotes}

\newsiamremark{remark}{Remark}
\newsiamremark{hypothesis}{Hypothesis}
\crefname{hypothesis}{Hypothesis}{Hypotheses}
\newsiamthm{claim}{Claim}

\headers{F.R.C. in V.R. Complexes and data applications}{D.B.S. et al}
\title{Efficient Decomposition of Forman-Ricci Curvature on Vietoris-Rips Complexes and Data Applications
\thanks{This manuscript is for review purposes only.
\funding{This research is supported by the grant\newline PID2023-146683OB-100 funded by MICIU/AEI /10.13039/501100011033 and by\newline ERDF, EU. Additionally, it is supported by Ikerbasque Foundation and the Basque Government through the BERC 2022-2025 program and by the Ministry of Science and Innovation: BCAM Severo Ochoa accreditation CEX2021-001142-S / MICIU / AEI / 10.13039/501100011033. Moreover, the authors acknowledge the financial support received from BCAM-IKUR, funded by the Basque Government by the IKUR Strategy and the European Union NextGenerationEU/PRTR. We also acknowledge the support of ONBODY no. KK-2023/00070 funded by the Basque Government through ELKARTEK Programme. Weiqiang Sun and Mengjun Ding are supported by the National Key Research and Development Project of China under Grant 2024YFB2908301 and by the National Natural Science Foundation of China (NSFC) under Grant 62331017}}}
\author{Danillo Barros de Souza\thanks{Basque Center for Applied Mathematics
  (\email{danillo.dbs16@gmail.com})}
  \and Jonatas Teodomiro\thanks{Universidade Federal de Pernambuco (\email{jonatas.teodomiro@ufpe.br})}
  \and Fernando A. N. Santos\thanks{Dutch Institute for Emergent Phenomenal (\email{f.a.nobregasantos@uva.nl})}
  \and
  Mengjun Ding\thanks{School of
  Electronic Information and Electrical Engineering,
  Shanghai Jiao Tong University, Shanghai, China (\email{mengjun\_ding@sjtu.edu.cn})}
  \and
  Weiqiang Sun\thanks{School of
  Electronic Information and Electrical Engineering,
  Shanghai Jiao Tong University, Shanghai, China (\email{sunwq@sjtu.edu.cn})}
  \and
  Mathieu Desroches\thanks{Inria, Montpellier, France (\email{mathieu.desroches@inria.fr
  })}
  \and J\"urgen Jost\thanks{Max Planck Institute for Mathematics in the Sciences, Leipzig,  and Center for Scalable Data Analytics and Artificial Intelligence, Leipzig University, Germany,  and Santa Fe Institute, New Mexico, USA
  (\email{jost@mis.mpg.de}).}
  \and Serafim Rodrigues\thanks{Basque Center for Applied Mathematics (\email{srodrigues@bcamath.org})}}

\usepackage{amsopn}
\usepackage{cleveref}

\makeatletter
\newcommand*{\addFileDependency}[1]{% argument=file name and extension
  \typeout{(#1)}% latexmk will find this if $recorder=0 (however, in that case, it will ignore #1 if it is a .aux or .pdf file etc and it exists! if it doesn't exist, it will appear in the list of dependents regardless)
  \@addtofilelist{#1}% if you want it to appear in \listfiles, not really necessary and latexmk doesn't use this
  \IfFileExists{#1}{}{\typeout{No file #1.}}% latexmk will find this message if #1 doesn't exist (yet)
}
\makeatother

%% file: references.bib
@article{fortunato2010community,
  title={Community detection in graphs},
  author={Fortunato, Santo},
  journal={Physics Reports},
  volume={486},
  number={3-5},
  pages={75--174},
  year={2010},
  publisher={Elsevier}
}

@article{newman2003structure,
  title={The structure and function of complex networks},
  author={Newman, Mark E. J.},
  journal={SIAM Review},
  volume={45},
  number={2},
  pages={167--256},
  year={2003},
  publisher={}}

@article{zhou2020graph,
  title={Graph neural networks: A review of methods and applications},
  author={Zhou, Jie and Cui, Ganqu and Hu, Shengding and Zhang, Zhengyan and Yang, Cheng and Liu, Zhiyuan and Wang, Lifeng and Li, Changcheng and Sun, Maosong},
  journal={AI open},
  volume={1},
  pages={57--81},
  year={2020},
  publisher={Elsevier}
}

@article{jin2021survey,
  title={A survey of community detection approaches: From statistical modeling to deep learning},
  author={Jin, Di and Yu, Zhizhi and Jiao, Pengfei and Pan, Shirui and He, Dongxiao and Wu, Jia and Philip, S Yu and Zhang, Weixiong},
  journal={IEEE Transactions on Knowledge and Data Engineering},
  volume={35},
  number={2},
  pages={1149--1170},
  year={2021},
  publisher={IEEE}
}

@article{zomorodian2012topological,
  title={Topological data analysis},
  author={Zomorodian, Afra},
  journal={Advances in Applied and Computational Topology},
  volume={70},
  pages={1--39},
  year={2012}
}

@article{gunnarcarlsson2018towards,
  title={Towards a new approach to reveal dynamical organization of the brain using topological data analysis},
  author={Saggar, Manish and Sporns, Olaf and Gonzalez-Castillo, Javier and Bandettini, Peter A and Carlsson, Gunnar and Glover, Gary and Reiss, Allan L},
  journal={Nature Communications},
  volume={9},
  number={1},
  pages={1--14},
  year={2018},
  publisher={Nature Publishing Group}
}

@article{edelsbrunner2008persistent,
  title={Persistent {H}omology -- a {S}urvey},
  author={Edelsbrunner, Herbert and Harer, John and others},
  journal={Contemporary Mathematics},
  volume={453},
  number={26},
  pages={257--282},
  year={2008},
  publisher={Providence, RI: American Mathematical Society}
}

@book{zomorodian2005topology,
  title={Topology for Computing},
  author={Zomorodian, Afra J},
  volume={16},
  year={2005},
  publisher={Cambridge University Press}
}

@book{enderton1977elements,
  title={Elements of set theory},
  author={Enderton, Herbert B},
  year={1977},
  publisher={Academic press}
}

@book{jech2003set,
  title={Set theory: The third millennium edition, revised and expanded},
  author={Jech, Thomas},
  year={2003},
  publisher={Springer}
}

@inproceedings{zimmermann1985applications,
  title={Applications of Fuzzy Set Theory to Mathematical Programming},
  author={Zimmermann, H.-J.},
  booktitle={Fuzzy Sets for Intelligent Systems},
  editor={Dubois, D. and Prade, H. and Yager, R.R.},
  pages={795--809},
  year={1993},
  publisher={Elsevier}
}

@article{pawlak2002rough,
  title={Rough set theory and its applications},
  author={Pawlak, Zdzis{\l}aw},
  journal={Journal of Telecommunications and Information Technology},
  pages={7--10},
  year={2002}
}

@unpublished{de2023efficient,
  title={An efficient set-theoretic algorithm for high-order {F}orman-{R}icci curvature},
  author={de Souza, D.B. and da Cunha, J.T.S. and Santos, F.A.N. and Jost, J. and Rodrigues, S.},
  note={{\it Proceedings of the Royal Society A} (in press)},
  year={2025}
}

@misc{desouza2025alternativesettheoreticalalgorithmsefficient,
      title={Alternative set-theoretical algorithms for efficient computations of cliques in Vietoris-Rips complexes}, 
      author={Danillo Barros de Souza and Jonatas Teodomiro and Fernando A. N. Santos and Mathieu Desroches and Serafim Rodrigues},
      year={2025},
      eprint={2502.14593},
      archivePrefix={arXiv},
      primaryClass={math.CO},
      url={https://arxiv.org/abs/2502.14593}, 
}

@article{forman2003bochner,
  title={Bochner's method for cell complexes and combinatorial {R}icci curvature},
  author={Forman, Robin},
  journal={Discrete and Computational Geometry},
  volume={29},
  number={3},
  pages={323--374},
  year={2003},
  publisher={Springer}
}

@article{Lin2011RicciGraphs,
    title = {Ricci curvature of graphs},
    year = {2011},
    journal = {Tohoku Mathematical Journal},
    author = {Lin, Yong and Lu, Linyuan and Yau, Shing-Tung},
    number = {4},
    pages = {605--627},
    volume = {63},
    url = {http://projecteuclid.org/euclid.tmj/1325886283},
    doi = {10.2748/tmj/1325886283},
    issn = {0040-8735}
}

@article{sandhu2016ricci,
  title={Ricci curvature: An economic indicator for market fragility and systemic risk},
  author={Sandhu, Romeil S and Georgiou, Tryphon T and Tannenbaum, Allen R},
  journal={Science Advances},
  volume={2},
  number={5},
  pages={e1501495},
  year={2016},
  publisher={American Association for the Advancement of Science}
}

@article{samal2018comparative,
  title={Comparative analysis of two discretizations of {R}icci curvature for complex networks},
  author={Samal, Areejit and Sreejith, RP and Gu, Jiao and Liu, Shiping and Saucan, Emil and Jost, J{\"u}rgen},
  journal={Scientific Reports},
  volume={8},
  number={1},
  pages={1--16},
  year={2018},
  publisher={Nature Publishing Group}
}

@article{chatterjee2021detecting,
  title={Detecting network anomalies using {F}orman--{R}icci curvature and a case study for human brain networks},
  author={Chatterjee, Tanima and Albert, R{\'e}ka and Thapliyal, Stuti and Azarhooshang, Nazanin and DasGupta, Bhaskar},
  journal={Scientific Reports},
  volume={11},
  number={1},
  pages={8121},
  year={2021},
  publisher={Nature Publishing Group UK London}
}

@article{de2021using,
  title={Using discrete {R}icci curvatures to infer {COVID}-19 epidemic network fragility and systemic risk},
  author={de Souza, Danillo Barros and Da Cunha, Jonatas TS and dos Santos, Everlon Figueir{\^o}a and Correia, Jailson B and da Silva, Hernande P and de Lima Filho, Jos{\'e} Luiz and Albuquerque, Jones and Santos, Fernando AN},
  journal={Journal of Statistical Mechanics: Theory and Experiment},
  volume={2021},
  number={5},
  pages={053501},
  year={2021},
  publisher={IOP Publishing}
}

@article{zomorodian2010fast,
  title={Fast construction of the {V}ietoris--{R}ips complex},
  author={Zomorodian, Afra},
  journal={Computers \& Graphics},
  volume={34},
  number={3},
  pages={263--271},
  year={2010},
  publisher={Elsevier}
}

@book{edelsbrunner2022computational,
  title={Computational Topology: An introduction},
  author={Edelsbrunner, Herbert and Harer, John L},
  year={2022},
  publisher={American Mathematical Society}
}

@article{knill2012discrete,
  title={A discrete Gauss-Bonnet type theorem},
  author={Knill, Oliver},
  journal={Elemente der Mathematik},
  volume={67},
  number={1},
  pages={1--17},
  year={2012}
}

@unpublished{knill2024gauss,
  title={Gauss-{B}onnet for {F}orm {C}urvatures},
  author={Knill, Oliver},
  note={{\it arXiv e-print} 2409.01425, \url{https://arxiv.org/abs/2409.01425}},
  year={2024}
}

@book{penrose2003random,
  title={Random Geometric Graphs},
  author={Penrose, Mathew},
  volume={5},
  year={2003},
  publisher={Oxford University Press}
}

@article{dall2002random,
  title={Random geometric graphs},
  author={Dall, Jesper and Christensen, Michael},
  journal={Physical Review E},
  volume={66},
  number={1},
  pages={016121},
  year={2002},
  publisher={APS}
}

@inproceedings{matejka2017same,
  title={Same stats, different graphs: generating datasets with varied appearance and identical statistics through simulated annealing},
  author={Matejka, Justin and Fitzmaurice, George},
  booktitle={Proceedings of the 2017 CHI conference on human factors in computing systems},
  pages={1290--1294},
  year={2017}
}

@Manual{datasaurus,
  title = {datasau{R}us: {D}atasets from the {D}atasaurus {D}ozen},
  author = {Colin Gillespie and Steph Locke and Rhian Davies and Lucy {D'Agostino McGowan}},
  year = {2024},
  note = {R package version 0.1.8, https://jumpingrivers.github.io/datasauRus/},
  url = {https://github.com/jumpingrivers/datasauRus},
}

@unpublished{dlotko2023persistence,
  title={Persistence {N}orms and the {D}atasaurus},
  author={Dlotko, Pawel and Rudkin, Simon},
  note={{\it arXiv e-print} 2309.13479, \url{https://arxiv.org/abs/2309.13479}},
  year={2023}
}

@misc{python,
  title = {Python {L}anguage {R}eference, version 3.13},
  author = {{Python Software Foundation}},
  year = {2023},
  note = {\url{https://docs.python.org/3/}},
}

@misc{networkx,
  author = {Aric A. Hagberg and Daniel A. Schult and Pieter J. Swart},
  title = {Network{X}: {N}etwork {A}nalysis in {P}ython},
  year = {2008},
  note = {\url{https://networkx.org}},
  howpublished = {Accessed: 2023-11-13}
}

@unpublished{umap,
  author = {McInnes, Leland and Healy, John and Melville, James},
  title = {{UMAP}: {U}niform {M}anifold {A}pproximation and {P}rojection for {D}imension {R}eduction},
  note = {{\it arXiv e-print} 1802.03426, \url{https://arxiv.org/abs/1802.03426}},
  year = {2018},
  url = {https://arxiv.org/abs/1802.03426}
}

@misc{umap-learn,
  author = {McInnes, Leland and Healy, John and Melville, James},
  title = {{UMAP}: {U}niform {M}anifold {A}pproximation and {P}rojection},
  year = {2018},
  note = {\url{https://umap-learn.readthedocs.io}}
}

@misc{gudhi,
  author = {The GUDHI Project},
  title = {{GUDHI}: {G}eometry {U}nderstanding in {H}igher {D}imensions},
  year = {2023},
  note = {Version 3.6.0, \url{https://gudhi.inria.fr}},
}

@article{gudhi:2014,
  author = {Jean-Daniel Boissonnat and Clément Maria},
  title = {The {GUDHI} {L}ibrary: {S}implicial {C}omplexes and {P}ersistent {H}omology},
  journal = {ACM Communications in Computer Algebra},
  volume = {47},
  number = {3/4},
  pages = {85--87},
  year = {2014}
}

@misc{METABRIC_Kaggle,
  author       = {Raghad Alharbi},
  title        = {Breast {C}ancer {G}ene {E}xpression {P}rofiles ({METABRIC})},
  year         = {2020},
  url          = {https://www.kaggle.com/datasets/raghadalharbi/breast-cancer-gene-expression-profiles-metabric},
  note         = {Accessed: 2024-11-13},
  publisher    = {Kaggle}
}

@misc{UCI_Wisconsin_Breast_Cancer,
  author       = {Dheeru, Dua and Karra Taniskidou, E.},
  title        = {Breast {C}ancer {W}isconsin ({D}iagnostic) {D}ata {S}et},
  year         = {2017},
  url          = {https://archive.ics.uci.edu/ml/datasets/breast+cancer+wisconsin+(diagnostic)},
  note         = {Accessed: 2024-11-13},
  publisher    = {UCI Machine Learning Repository}
}

@article{bochner1949curvature,
  title={Curvature and {B}etti numbers. {II}},
  author={Bochner, Salomon},
  journal={Annals of Mathematics},
  volume={50},
  number={1},
  pages={77--93},
  year={1949},
  publisher={JSTOR}
}

@unpublished{knill2020index,
  title={On index expectation curvature for manifolds},
  author={Knill, Oliver},
  note={{\it arXiv e-print} 2001.06925, \url{https://arxiv.org/abs/2001.06925}},
  year={2020}
}

@inproceedings{xu2022amcad,
  title={AMCAD: adaptive mixed-curvature representation based advertisement retrieval system},
  author={Xu, Zhirong and Wen, Shiyang and Wang, Junshan and Liu, Guojun and Wang, Liang and Yang, Zhi and Ding, Lei and Zhang, Yan and Zhang, Di and Xu, Jian and others},
  booktitle={2022 IEEE 38th International Conference on Data Engineering (ICDE)},
  pages={3439--3452},
  year={2022},
  organization={IEEE}
}

@misc{FRC_VR_kaggle,
  author       = {Danillo Barros de Souza},
  title        = {{F.R.C.} on {V.R.} filtrations of random graphs},
  year         = {2025},
  howpublished = {\url{https://kaggle.com/datasets/fa3926660ecbe1ced3e2de6012a32ba9e3bdfd8988f1e6ec9250387aafc73214}},
  note         = {Accessed: 2025-04-24}
}

@article{bloch2014combinatorial,
  title={Combinatorial Ricci curvature for polyhedral surfaces and posets},
  author={Bloch, Ethan},
  journal={arXiv preprint arXiv:1406.4598},
  year={2014}
}

@article{weber2018coarse,
  title={Coarse geometry of evolving networks},
  author={Weber, Melanie and Saucan, Emil and Jost, J{\"u}rgen},
  journal={Journal of complex networks},
  volume={6},
  number={5},
  pages={706--732},
  year={2018},
  publisher={Oxford University Press}
}

@misc{souza2025fastforman_benchmark,
  author       = {Danillo Barros de Souza},
  title        = {Example Usage and Benchmarks},
  howpublished = {\url{https://github.com/danillodbs16/fastforman/blob/main/codes/0.7.0/Benchmarks/Examples/Example_usage_and_Benchmarks.ipynb}},
  year         = {2025},
  note         = {Jupyter notebook, FastForman GitHub repository. Accessed: 2026-08-03}
}

@software{Barros_de_Souza_FastForman_-_An_2024,
author = {Barros de Souza, Danillo},
doi = {10.5281/zenodo.11167956},
month = may,
title = {{FastForman - An efficient Forman-Ricci Curvature computation for higher-order faces in Simplicial Complexes}},
url = {https://github.com/danillodbs16/fastforman},
version = {0.1.8},
year = {2024}
}

@article{jost2014ollivier,
  title={Ollivier’s Ricci curvature, local clustering and curvature-dimension inequalities on graphs},
  author={Jost, J{\"u}rgen and Liu, Shiping},
  journal={Discrete \& Computational Geometry},
  volume={51},
  number={2},
  pages={300--322},
  year={2014},
  publisher={Springer}
}

@inproceedings{shahapure2020cluster,
  title={Cluster quality analysis using silhouette score},
  author={Shahapure, Ketan Rajshekhar and Nicholas, Charles},
  booktitle={2020 IEEE 7th international conference on data science and advanced analytics (DSAA)},
  pages={747--748},
  year={2020},
  organization={IEEE}
}

@software{ni2024graphriccicurvature,
  author       = {Chien-Chun Ni},
  title        = {{GraphRicciCurvature}: A Python Library to Compute Graph Ricci Curvature and Ricci Flow},
  year         = {2024},
  version      = {0.5.3.2},
  url          = {https://github.com/saibalmars/GraphRicciCurvature},
  note         = {GitHub repository, accessed August 4, 2026}
}

@inproceedings{cohen2005stability,
  title={Stability of persistence diagrams},
  author={Cohen-Steiner, David and Edelsbrunner, Herbert and Harer, John},
  booktitle={Proceedings of the twenty-first annual symposium on Computational geometry},
  pages={263--271},
  year={2005}
}
